\documentclass[11pt]{amsart}
\usepackage{pst-all}
\usepackage{amsmath}
\usepackage[latin1]{inputenc}
\usepackage{amsfonts}
\usepackage{amssymb}

\usepackage{amsthm}
\usepackage{amscd}

\usepackage[latin1]{inputenc}

\usepackage{amsmath}

\usepackage{amsfonts}

\usepackage{amssymb}

\usepackage{amsthm}

\usepackage{graphics}

\newcommand{\mk}{\medskip}

\newcommand{\ZZ}{\mathbb{Z}}

\newcommand{\CC}{\mathbb{C}}

\newcommand{\Glie}{\mathfrak{g}}

\newcommand{\demo}{\noindent {\it \small Proof:}\quad}

\newcommand{\U}{\mathcal{U}}

\newcommand{\nc}{\newcommand}
\nc{\g}{{\mathfrak g}}
\nc{\ghat}{\widehat{\g}}
\nc{\mc}{\mathcal}
\nc{\ep}{\epsilon}
\nc{\la}{\lambda}
\nc{\Z}{{\mathbb Z}}
\nc{\C}{{\mathbb C}}
\nc{\on}{\operatorname}
\nc{\wt}{\widetilde}
\nc{\om}{\omega}
\nc{\ol}{\overline}

\newtheorem{thm}{Theorem}[section]

\newtheorem{defi}[thm]{Definition}

\newtheorem{prop}[thm]{Proposition}

\newtheorem{lem}[thm]{Lemma}

\newtheorem{conj}[thm]{Conjecture}

\newtheorem{rem}[thm]{Remark}

\title{Langlands duality for representations of quantum groups}

\author[Edward Frenkel]{Edward Frenkel$^1$}\thanks{$^1$Supported in
  part by DARPA and AFOSR through the grant FA9550-07-1-0543}

\address{Department of Mathematics, University of California,
  Berkeley, CA 94720, USA}

\author[David Hernandez]{David Hernandez$^2$}\thanks{$^2$Supported partially by ANR through Project "G\'eom\'etrie et Structures
  Alg\'ebriques Quantiques"}

\address{CNRS - 
\'Ecole Normale Sup\'erieure, DMA, 45, rue d'Ulm 75005 Paris,
  FRANCE}

\date{September 2008}

\usepackage[all]{xy}
\begin{document}

\begin{abstract} We establish a correspondence (or duality) between the
characters and the crystal bases of finite-dimensional representations
of quantum groups associated to Langlands dual semi-simple Lie
algebras. This duality may also be stated purely in terms of
semi-simple Lie algebras.  To explain this duality, we introduce an
``interpolating quantum group'' depending on two parameters which
interpolates between a quantum group and its Langlands dual. We
construct examples of its representations, depending on two
parameters, which interpolate between representations of two Langlands
dual quantum groups.

\vskip 4.5mm


\end{abstract}

\maketitle


\pagestyle{myheadings}

\section{Introduction}

Let $\g$ be a simple Lie algebra and $^L\g$ its Langlands dual Lie
algebra whose Cartan matrix is the transpose of that of $\g$. In this
paper we establish a {\em duality} between finite-dimensional
representations of $\g$ and $^L\g$, as well as the corresponding
quantum groups.

Let $I$ be the set of vertices of the Dynkin diagram of $\g$ and $r_i,
i \in I$, the corresponding labels. Denote by $r$ the maximal number
among the $r_i$. This is the lacing number of $\g$ which is equal to
$1$ for the simply-laced $\g$, to $2$ for $B_\ell, C_\ell$ and $F_4$,
and to $3$ for $G_2$.

Let $L(\la)$ be a finite-dimensional irreducible representation of
$\g$ whose highest weight $\la$ has the form
\begin{equation}    \label{la}
\la = \sum_{i \in I} (1+r-r_i) m_i \omega_i, \qquad m_i \in \Z_+,
\end{equation}
where the $\omega_i$ are the fundamental weights of $\g$. In other
words, $\la$ is a dominant integral weight which belongs to the
sublattice $P' \subset P$, where $P$ is the weight lattice of $\g$,
spanned by $(1+r-r_i) \omega_i, i \in I$. The character of $L(\la)$
has the form
$$
\chi(L(\la)) = \sum_{\nu \in P} d(\la,\nu) e^\nu, \qquad d(\la,\nu)
\in \Z_+.
$$
Let
$$
\chi'(L(\la)) = \sum_{\nu \in P'} d(\la,\nu) e^\nu.
$$
We first prove that, after replacing each
$$
\nu = \sum_{i \in I} (1+r-r_i) n_i \omega_i \, \in \, P', \qquad n_i
\in \Z,
$$
by
$$
\nu' = \sum_{i \in I} n_i \check\omega_i,
$$
where the $\check\omega_i$ are the fundamental weights of $^L\g$,
$\chi'(L(\la))$ becomes the character of a virtual finite-dimensional
representation of $^L\g$, whose highest component is $L(\la')$, the
irreducible representation of $^L\g$ with the highest weight
\begin{equation}    \label{lap}
\la' = \sum_{i \in I} m_i \check\omega_i,
\end{equation}
where the numbers $m_i$ are defined by formula \eqref{la}. In other
words, we have
\begin{equation}    \label{branching}
\chi'(L(\la)) = \chi^L(L(\lambda')) + \sum_{\check\mu<\la'}
m_{\check\mu} \; \chi^L(L(\check\mu)), \qquad m_{\check\mu} \in
\Z.
\end{equation}
Then we prove that the multiplicities of weights in the character
$\chi^L(L(\lambda'))$ of $L(\lambda')$ are less than or equal to those
in $\chi'(L(\la))$. This positivity result means that
$\chi^L(L(\lambda'))$ is ``contained as a subcharacter'' in
$\chi'(L(\lambda))$.

\medskip

Since the categories of finite-dimensional representations of $\g$ and
$U_q(\g)$ with generic $q$ are equivalent, we also obtain a
duality between finite-dimensional representations of $U_q(\g)$ and
$U_q({}^L\g)$. Moreover, we establish the duality not only at the
level of characters but at the level of {\em crystal bases} as well.
This leads, in particular, to the following surprising fact: one can
construct the crystal basis of the irreducible representation
$L(\la')$ of $^L \g$ from the crystal basis of the irreducible
representation $L(\la)$ of $\g$. 
\footnote{After the first version of this paper appeared on the arXiv,
we learned from Hiraku Nakajima that this result follows from a special
case of \cite[Theorem 5.1]{kas0}; see the paragraph before Theorem
\ref{dcrys} for more details.}

In addition, we conjecture that $\chi'(L(\la))$ is the character of an
{\em actual} representation of $^L\g$ (that is, $m_{\check\mu} \geq 0$
for all $\check\mu$ in formula \eqref{branching}), and we prove this
conjecture for $\g=B_2$.\footnote{After the first version of this
  paper appeared on the arXiv we learned from Victor Kac that a
  special case of our duality, going from type $B$ to type $C$, may be
  explained in the context of representation theory of Lie
  superalgebras of type $B(0,n)$ as defined in \cite{kac}. In fact,
  the condition on the highest weight $\lambda\in P'$ appears in this
  case in \cite[Theorem 8]{kac} in the form $a_n\in 2\ZZ$. It is not
  clear to us whether one can use Lie superalgebras to interpret our
  duality for other types. \\In addition, we have learned from Kevin
  McGerty that in the meantime he has been able to prove this
  conjecture for other types, see \cite{M}.  \\After the paper was
  published, we were informed by C\'edric Lecouvey that this result
  was obtained earlier by Peter Littelmann \cite[Proposition 4]{pl}.}
We observe that the subset of the crystal of $L(\la)$ consisting of
those elements whose weights are in $P'$ does not give us ``on the
nose'' the crystal of this $^L\Glie$-module. But we conjecture that
after applying a certain deformation process (presented in Section
\ref{conj}) we do get the right crystal structure on this subset. (We
also prove this for $B_2$.)  Thus, conjecturally, we can reconstruct
not only the crystal of $L(\la')$ but the crystal of the whole
representation of $^L\g$ whose character is equal to $\chi'(L(\la))$.

\medskip

It is natural to ask: {\em why} should this duality of characters and
crystal bases take place? We suggest the following explanation: there
exists an algebra $\U_{q,t}(\g)$ depending on two parameters, $q$ and
$t$, whose specialization at $t=1$ gives $\U_q(\g)$, and at $q=\ep$
(where $\ep=1$ if $\g$ is simply-laced and $\ep = \exp(\pi i/r)$, $r$
being the lacing number of $\g$) gives $\U_{-t}({}^L\g)$. These are
the quantum groups without the Serre relations associated to $\g$ and
$^L\g$. We call $\U_{q,t}(\g)$ the {\em interpolating quantum
group}. (Example 3 in Section \ref{int reps} indicates that it is
impossible to include the Serre relations and preserve the
interpolating property.) Moreover, we conjecture that any irreducible
finite-dimensional representation $L_q(\la)$ of $U_q(\g)$
(equivalently, of $\U_q(\g)$) with the highest weight of the form
\eqref{la} may be deformed to a representation $L_{q,t}(\la)$ of
$\U_{q,t}(\g)$. We also conjecture that the specialization of
$L_{q,t}(\la)$ at $q=\ep$ contains the irreducible representation of
$U_{-t}({}^L\g)$ with highest weight $\la'$ given by formula
\eqref{lap} as the highest component. These conjectures are confirmed
by various explicit examples presented below as well as our general
result on the duality of characters of finite-dimensional
representations.

\bigskip

Now we would like to briefly sketch a possible link between our
results and the geometric Langlands correspondence (see, e.g,
\cite{F:rev} for a general introduction).

One of the key results used in the geometric Langlands correspondence
is an isomorphism between the center $Z(\ghat)$ of the completed
enveloping algebra of $\ghat$ at the critical level and the classical
${\mc W}$-algebra ${\mc W}({}^L\g)$ (see \cite{FF:gd,F:wak} as well as
\cite{fr3} for details). This result forms the basis for the local
geometric Langlands correspondence (see \cite{FG,f}) as well as for
the Beilinson--Drinfeld construction of the global geometric Langlands
correspondence \cite{BD} (see also \cite{F:rev}). However, this
isomorphism is rather mysterious. We know that it exists but we do
not fully understand {\em why} it should exist.

In order to understand this better, we $q$-deform the picture and
consider the center $Z_q(\ghat)$ of the quantum affine algebra
$U_q(\ghat)$ at the critical level, which was the starting point of
\cite{fr}. The center $Z_q(\ghat)$ is in turn related to the
Grothendieck ring $\on{Rep} U_q(\ghat)$ of finite-dimensional
representations of $U_q(\ghat)$ (this is because for each
finite-dimensional representation $V$ we can construct a generating
series of central elements in $Z_q(\ghat)$, using the transfer-matrix
construction). Thus, we hope to gain some insight into the isomorphism
$Z(\ghat) \simeq {\mc W}({}^L\g)$ by analyzing the connections between
$Z_q(\ghat)$, $\on{Rep} U_q(\ghat)$ and the $q$-deformed classical
${\mc W}$-algebra.

The idea of \cite{fr2} was to further deform this picture and
introduce a two-parameter (non-commutative) deformation ${\mc
W}_{q,t}(\g)$. Its specialization ${\mc W}_{q,1}(\g)$ at $t=1$ is the
center $Z_q(\ghat)$, so that ${\mc W}_{q,t}(\g)$ is a one-parameter
deformation of $Z_q(\ghat)$ and a two-parameter deformation of the
original center $Z(\ghat)$. The work \cite{fr2} was motivated by the
hope that analyzing various dualities and limits of ${\mc
W}_{q,t}(\g)$ we may learn something new about the isomorphism
$Z(\ghat) \simeq {\mc W}({}^L\g)$ and hence about the Langlands
correspondence.

In particular, it was suggested in \cite{fr2} that the specialization
${\mc W}_{\ep,t}(\g)$ at $q=\ep$ (with $\ep$ defined as above)
contains as a subalgebra the center $Z_t({}^L\ghat)$ of the quantum
affine algebra $U_t({}^L\ghat)$ at the critical level (here $^L\ghat$
denotes the Langlands dual of $\ghat$). The latter gives rise to the
Grothendieck ring of finite-dimensional representations of
$U_t({}^L\ghat)$ (via the transfer-matrix construction). On the other
hand, as we already mentioned above, the specialization ${\mc
W}_{q,1}(\g)$ at $t=1$ gives rise to the Grothendieck ring of
finite-dimensional representations of $U_q(\ghat)$. Thus, it appears
that the ${\mc W}$-algebra ${\mc W}_{q,t}(\g)$ {\em interpolates}
between the Grothendieck rings of finite-dimensional representations
of quantum affine algebras associated to $\ghat$ and $^L\ghat$. In
particular, this suggests that these representations should be related
in some way. Examples of such a relation were given in \cite{fr2}, but
this phenomenon has largely remained a mystery until now.

How can we explain this relation from the point of view of
representation theory? This question served as the motivation for this
paper. Before answering it, we considered its {\em finite-dimensional
analogue}: is there a hidden correspondence, or duality, between
finite-dimensional representations of the quantum groups $U_q(\g)$ and
$U_q({}^L \g)$ -- or the simple Lie algebras $\g$ and $^L\g$, for that
matter?

We have given an affirmative answer to this question which we have
outlined above. Thus, we have found a hidden duality between objects of
the same nature: finite-dimen\-sional representations of two Langlands
dual Lie algebras. Actually, it is rather surprising to observe the
appearance of a Langlands type duality in such an elementary context:
that of finite-dimensional representations of simple Lie algebras! We
hope that this duality and its affine analogue will give us some clues
about the meaning of the geometric Langlands correspondence.

\medskip

What about the duality for the quantum affine algebras? In our next
paper \cite{FH:new} we will propose a precise relation between the
$q$-characters of finite-dimensional representations of dual quantum
affine algebras $U_q(\ghat)$ and $U_q(^L\ghat)$ with is analogous to
the duality of characters of $U_q(\g)$ and $U_q({}^L\g)$ discussed
above. We will prove, by using \cite{hcr, hgen}, that this relation
holds for an important class of representations, the {\em
Kirillov--Reshetikhin modules}. In the affine case we also expect that
the duality may be explained by using an affine analogue of the
interpolating quantum group.

\medskip

In the context of our results an interesting problem is to compute
explicitly all multiplicities of simple $^L\g$-modules in a given
simple $\g$-module (the numbers $m_{\check\mu}$ in formula
\eqref{branching}), which we call the {\em Langlands duality branching
rules}.  In the course of the proof we have found them explicitly in
some cases.

\medskip

The paper is organized as follows. In Section \ref{duality of reps} we
establish the duality of characters and crystal bases for a pair of
Langlands dual simple Lie algebras. In Section \ref{iqg} we introduce
the interpolating quantum group. We then study its representations
which we expect to interpolate between representations of $U_q(\g)$
and $U_{-t}({}^L\g)$. This would explain the duality that we have
found in this paper. In Section \ref{bun} we show how this
interpolation works for the finite-dimensional representations of the
elementary interpolating quantum groups (those corresponding to Lie
algebras of rank one). In Section \ref{int reps} we consider examples
of more general interpolating representations. In Section \ref{conj}
we conjecture a stronger duality for characters and crystals and
prove it for all simply-laced $\g$ with $r = 2$ and for $B_2$.

\medskip

\noindent{\bf Acknowledgments.} This work was begun while we were
taking part in the Program on Combinatorial Representation Theory held
at MSRI in the Spring of 2008. We thank the organizers of this Program
for their invitations and MSRI for hospitality.

\section{Duality of characters and crystals for simple Lie algebras}
\label{duality of reps}

In this section we prove the Langlands duality for characters of
finite-dimensional representations of quantum groups associated to
simple Lie algebras (or, equivalently, simple Lie algebras
themselves). We also prove the duality of the corresponding crystal
bases, by using the monomial model \cite{Nac, kas}.

\medskip

Let $\g$ be a finite-dimensional simple Lie algebra and $U_q(\g)$ the
corresponding quantum group (see, e.g., \cite{cp}). We denote $r =
\text{max}_{i\in I}(r_i)$, where $I$ is the set of vertices of the
Dynkin diagram of $\g$ and the $r_i$ are the corresponding
labels. This is the lacing number of $\g$ (note that it was denoted by
$r^\vee$ in \cite{fr2}).

The Cartan matrix of $\g$ will be denoted by $C=(C_{i,j})_{i,j \in
I}$. By definition, the Langlands dual Lie algebra $^L\g$ has the
Cartan matrix $C^t$, the transpose of the Cartan matrix $C$ of $\g$.

\subsection{Langlands duality for characters}\label{charstate}

Let $$P=\sum_{i \in I} \ZZ\om_i$$ be the weight lattice of $\Glie$ and
$P^+ \subset P$ the set of dominant weights. Consider the sublattice
\begin{equation}    \label{Pprime}
P' = \sum_{i\in I} (1 + r - r_i)\ZZ\om_i \subset P.
\end{equation}
Let $$P^L=\sum_{i \in I} \ZZ\check\om_i$$ be the weight lattice of
$^L\g$. Consider the map $\Pi : P\rightarrow P^L$ defined by
$$\Pi(\lambda) = \sum_{i\in I}\lambda(\check\alpha_i)(1 + r -
r_i)^{-1} \check\omega_i$$ if $\lambda \in P'$ and $\Pi(\lambda) = 0$,
otherwise. Clearly, $\Pi$ is surjective.

In this section we investigate what $\Pi$ does to characters of
irreducible representations of $\g$. For simply-laced Lie algebras
(that is, $r=1$) we have $P' = P = P^L$ and $\Pi$ is the
identity. Hence we focus on the non-simply laced Lie algebras.

Let $\on{Rep} \g$ be the Grothendieck ring of finite-dimensional
representations of $\g$. We have the character homomorphism
$$
\chi: \on{Rep} \g \to \Z[P] = \Z[y_i^{\pm 1}],
$$
where $y_i = e^{\om_i}$. It sends an irreducible representation
$L(\la)$ of $\g$ with highest weight $\la \in P^+$ to its character,
which we will denote by $\chi(\la)$.

We will now show that for any representation $V$ of
$\g$, $\Pi(\chi(V))$ is the character of a virtual
representation of $^L\g$, as stated in the following proposition. 
We denote the character homomorphism for $^L\g$ by $\chi^L$.

\begin{prop}\label{debut} For any simple Lie algebra $\g$ and any
  $\lambda \in P^+$, $\Pi(\chi(\lambda))$ is in the image of $\chi^L$.
\end{prop}

This is a direct consequence of the following Lemma. Here we denote by
$s_i$ (resp. $s_i^L$) the simple reflections of $\Glie$
(resp. $^L\g$).

\begin{lem} $P'$ is invariant under the Weyl group action and $\Pi
  \circ s_i = s_i^L \circ \Pi$ on $P'$.
\end{lem}

\demo Let $\mu = \prod_{j\in I}y_j^{\mu_j}\in P'$ and $i\in I$.

If $r_i = r$, we have $s_i(\mu) = \mu y_i^{-2\mu_i}(\prod_{j\sim i,
  r_j = 1}y_j^{r\mu_i})(\prod_{j\sim i, r_j = r} y_i^{\mu_i})\in
P'$. Moreover $\Pi(s_i(\mu)) = \Pi(\mu) y_i^{-2\mu_i}(\prod_{j\sim
  i}y_j^{\mu_i}) = s_i^L(\Pi(\mu))$.

If $r_i = 1$, we have $\mu_i\in r\ZZ$ and $s_i(\mu) = \mu
y_i^{-2\mu_i}(\prod_{j\sim i}y_j^{\mu_i})\in P'$. Moreover
$\Pi(s_i(\mu)) = \Pi(\mu) y_i^{-2\mu_i}(\prod_{j\sim i,r_j =
  1}y_j^{\mu_i/r})(\prod_{j\sim i,r_j = r}y_j^{\mu_i}) =
s_i^L(\Pi(\mu))$.
\qed

\begin{rem} If $\Glie$ is of type $B_\ell$ and $\lambda\in P^+\cap
  P'$, then all terms in $\chi(\la)$ correspond to weights in
$P'$, and so $\Pi(\chi(\la))$ has the same number of monomials
as $\chi(\la)$.
\end{rem} 

According to Proposition \ref{debut}, we have, for $\la \in P' \cap
P^+$,
$$
\Pi(\chi(\lambda)) = \sum_{\check\mu \in P^{L,+}} m_{\check\mu}
\chi^L(\check\mu), \qquad m_{\check\mu} \in \Z.
$$
It is clear from the definition that the maximal $\check\mu$ for which
$m_{\check\mu} \neq 0$ is the image of $\la$ under $\Pi$. Moreover, in
this case $m_{\check\mu}  = 1$. An interesting problem is
to compute explicitly all other multiplicities $m_{\check\mu}$, the
{\em Langlands duality branching rules}.

One of the main results of this section is the following:

\begin{thm}\label{incchar}
The multiplicities of weights in $\chi^L(\Pi(\lambda))$ are less than
or equal to those in $\Pi(\chi(\lambda))$.
\end{thm}

In other words, $\chi^L(\Pi(\lambda))$ can be seen as a
``subcharacter'' contained in $\Pi(\chi(\lambda))$, that is,
$\chi^L(\Pi(\lambda))\preceq \Pi(\chi(\lambda))$ where $\preceq$ is
the obvious partial ordering on polynomials.

\begin{rem} In general, the character $\chi(\lambda)$ is given by
  the Weyl character formula. So one could try to prove the above
  results by using the Weyl formula. However, it is not clear how to
  do this: although the Weyl groups of $\Glie$ and $^L\g$ are
  isomorphic, there is no obvious relation for the half-sums of
  positive roots $\rho$ and $\check\rho$.
\end{rem}

Before giving the proof, we consider some explicit examples.

Let $\Glie = B_2$. Then $^L\g = C_2$, which is isomorphic to $B_2$ but
with the switch of the labels of the Dynkin diagram $1 \rightarrow
\ol{1} = 2, 2 \rightarrow \ol{2}=1$. In other words, $\check\om_i$
corresponds not to $\omega_i$, but to $\om_{\ol{i}}$.

We have $P' = \ZZ\om_1 + 2\ZZ\om_2$. Here are the simplest examples of
action of $\Pi$ on characters of irreducible representations:

$$\Pi(\chi(\om_1)) = (y_1 + y_2y_1^{-1} + y_1 y_2^{-1} + y_1^{-1}) +
1 \succeq \chi^L(\check\om_1).$$

$$\Pi(\chi(2\om_2)) = (y_2 + y_2^{-1}y_1^2 + 1 + y_1^{-2}y_2 +
y_2^{-1}) + y_1 + y_2y_1^{-1} + 1 + y_2^{-1}y_1 + y_1^{-1}  \succeq
\chi^L(\check\om_2).$$

$$\Pi(\chi(2\om_1)) = (y_1^2 + y_2 + y_1^{-2}y_2 + y_1^2 y_2^{-1} +
2 + y_1^{-2}y_2 + y_2^{-2}y_1^2 + y_2^{-1} + y_1^{-2}) $$ $$+ y_1 +
y_1^{-1} y_2 + y_2^{-1}y_1 + y_1^{-1} \succeq \chi^L(2\check\om_1).$$

Let us look at some examples for $\Glie=G_2$. In this
case $^L\g=G_2$, but again with the switch of labels of the Dynkin
diagram, as in the case of $B_2$.
We have $P' = \ZZ\om_1 + 3\ZZ\om_3$. Here are a few examples:
$$\Pi(\chi(\om_1)) = (y_1 + y_2 y_1^{-1} + y_2^{-1}y_1^2 + 1 +
y_1^{-2} y_2 + y_2^{-1}y_1 + y_1^{-1}) + 1
\succeq \chi^L(\check\om_1).$$
$L(3\om_2)$ is of dimension $77$. We will not write it out explicitly,
but only write
$$\Pi(\chi(\om_2)) = 
(y_2
+ y_2^{-1}y_1^3 
+ y_1  
+ y_1^{-1}y_2  
+ y_2^2y_1^{-3}  
+ y_1^2y_2^{-1} 
$$ $$+ 2
+ y_1^{-2}y_2
+ y_1^3y_2^{-2}
+ y_1y_2^{-1}
+ y_1^{-1}
+ y_2y_1^{-3}
+ y_2^{-1})
$$ $$+
2 y_1  
+ 2 y_1^{-1}y_2  
+ 2 y_1^2y_2^{-1} 
+ 3
+ 2 y_1^{-2}y_2
+ 2 y_1y_2^{-1}
+ 2 y_1^{-1}
\succeq \chi^L(\check\om_2).$$

\subsection{Langlands duality of crystals of irreducible
  representations}

To prove Theorem \ref{incchar}, we will use the crystal basis theory. It
gives us an algorithm to compute character formulas. We will see that
the statement of Theorem \ref{incchar} is actually satisfied at the
level of crystal. Before proving this, we state a closely related
result describing a duality of crystals of irreducible representations
of $U_q(\g)$ and $U_q({}^L\g)$.

Let $\lambda\in P'\cap P^+$ and $\mathcal{B}(\lambda)$ be the
corresponding crystal of $L(\la)$, with a highest element
$u_\lambda$ and crystal operators $e_i, f_i$. We consider the operators
\begin{equation}    \label{tilde f e}
f_i^L = f_i^{1 + r - r_i}\text{ , }e_i^L =
e_i^{1 + r - r_i}.
\end{equation}
Let $\mathcal{B}'(\lambda)$ be the connected component of $u_\lambda$
in $\mathcal{B}(\lambda)$ for the operators $f_i^L,
e_i^L$. Note that the definition of $\mathcal{B}'(\lambda)$
depends only on the structure of the $\Glie$-crystal of
$\mathcal{B}(\lambda)$. The weight of the elements of
$\mathcal{B}'(\lambda)$ are in $P'$ and so for
$v\in\mathcal{B}'(\lambda)$ we can define $\mathrm{wt}^L(v) =
\Pi(\mathrm{wt}(v))$. Then for any simple Lie algebra $\g$ (including
$G_2$) we have the following theorem.

After the first version of this paper appeared, we learned from
H. Nakajima that this theorem follows from a special case of
\cite[Theorem 5.1]{kas0} (namely, we put $\xi = \text{Id}_I$ and $m_i
= 1 + r - r_i$ in the notation of \cite{kas0}). Note that \cite{kas0}
discussed examples of embeddings $\mathcal{B}(\lambda)\rightarrow
\mathcal{B}(m\lambda)$ and foldings obtained from automorphisms of
simply-laced Dynkin diagrams, whereas in the present paper we view
this in the context of Langlands duality.

\begin{thm}\label{dcrys} For $\lambda\in P'\cap P^+$,
  $(\mathcal{B}'(\lambda),e_i^L,f_i^L,\mathrm{wt}^L)$
  is isomorphic to the $^L\g$-crystal
  $\mathcal{B}^L(\Pi(\lambda))$ of $L(\Pi(\la))$.
\end{thm}

Thus, by using only the crystal of the
$\Glie$-module $L(\la)$ we have constructed the crystal of the
$^L\g$-module $L(\Pi(\la))$. 

\begin{rem} Let us look at $\g=B_2$. If $p$ is even, to the
  representation $L(m\om_1 + p\om_2)$ of $\Glie$ corresponds the
  representation $L(m\om_1 + p\om_2/2)$ of $^L\g=C_2$. But
  $C_2 \simeq B_2$. So if in addition $m$ is even,
  to the representation $L(m\om_1 + p\om_2/2)$ of $^L\g$ corresponds
  the representation $L(m\om_1/2 + p\om_2/2)$ of ${}^L({}^L\g) =
  \Glie$. Thus, we see that this Langlands duality here is not an
  involution.
\end{rem}

Theorem \ref{dcrys} implies Theorem \ref{incchar} as we have
$$\Pi\left(\sum_{m'\in\mathcal{B}'(\lambda)} \mathrm{wt}(m)\right) =
\chi^L(\Pi(\lambda)).$$ 

\subsection{Reminder -- monomial crystals}

Let $C$ be a Cartan matrix of finite type and $s\colon I\rightarrow
\{0,1\}$ ($i\mapsto s_i$) a map such that $C_{i,j}\leq -1$ implies
$s_i+s_j=1$. Introduce formal variables $Y_{i,l}$, and let $A$ be the
set of monomials of the form $$m={\prod_{i\in
I,l\in\ZZ}}Y_{i,l}^{u_{i,l}(m)}, \qquad u_{i,l}(m)\in\ZZ.$$
A monomial $m$ is said to be dominant if $\forall j\in I, l\in\ZZ, u_{j,l}(m)\geq 0$. We set
$$A_{i,l}=Y_{i,l-1}Y_{i,l+1}{\prod_{j\neq i}}Y_{j,l}^{C_{j,i}}\in A.$$

Consider the subgroup $\mathcal{M}\subset A$ defined by
$$\mathcal{M}=\{m\in A\mid u_{i,l}(m)=0\text{ if }l\equiv s_i + 1
\text{ mod } 2\}.$$

Let us define $\mathrm{wt}\colon A\to P$ and $\epsilon_i,\phi_i,
p_i,q_i\colon A\to \ZZ\cup \{\infty\}\cup\{-\infty\}$, 
$e_i,f_i\colon A\to A\cup\{0\}$ for $i\in I$ by
the formulas (for $m\in A$)
\begin{gather*}
\mathrm{wt}(m)=\sum_{i\in I, l\in\ZZ}u_{i,l}(m) \omega_i,
\\
\phi_{i,L}(m) = \sum_{l\leq L} u_{i,l}(m), 
\quad
\phi_i(m)=\max\{0,\{ \phi_{i,L}(m) \mid L\in\ZZ\}\}\geq 0,
\\
\epsilon_{i,L}(m)=-{\sum_{l\geq L}}u_{i,l}(m),
\quad
\epsilon_i(m) = \max\{ 0,\{\epsilon_{i,L}(m) \mid L\in\ZZ\}\}\geq 0,
\\
p_i(m)=\max\{L\in\ZZ\mid
\epsilon_{i,L}(m)=\epsilon_i(m)\}=\max\{L\in\ZZ\mid
	{\sum_{l<L}}u_{i,l}(m)=\phi_i(m)\},
\\
q_i(m)=\text{min}\{L\in\ZZ\mid
\phi_{i,L}(m)=\phi_i(m)\}=\text{min}\{L\in\ZZ\mid -{\sum_{l>
    L}}u_{i,l}(m)=\epsilon_i(m)\}.
\end{gather*}
\begin{equation*}
  e_i(m) =
  \begin{cases}
    0 & \text{if $\epsilon_i(m) = 0$},
    \\
   mA_{i,p_i(m)-1} & \text{if $\epsilon_i(m) > 0$},
  \end{cases}
\text{ and } f_i(m) =
  \begin{cases}
    0 & \text{if $\phi_i(m) = 0$},
    \\
   mA_{i,q_i(m)+1}^{-1} & \text{if $\phi_i(m) > 0$}.
  \end{cases}
\end{equation*}
By \cite{Nac,kas} $(\mathcal{M},\mathrm{wt}, \epsilon_i, \phi_i,
e_i, f_i)$ is a crystal (called the monomial
crystal). For $m\in \mathcal{M}$ we denote by $\mathcal{M}(m)$ the
subcrystal of $\mathcal{M}$ generated by $m$.

\begin{thm}\label{real}\cite{Nac,kas} If $m$ is dominant, then the
  crystal $\mathcal{M}(m)$ is isomorphic to the crystal
  $\mathcal{B}(\mathrm{wt}(m))$ of $L(\mathrm{wt}(m))$.\end{thm}

In the following we will use the notation $i_l^r$ for $Y_{i,l}^r$.

\subsection{Examples}\label{excrys} We first study examples for Lie
algebras of rank $2$ and the following representations:

\begin{defi} The irreducible representations $L((r + 1 - r_i)\om_i)$
  will be called {\em pseudo fundamental representations}, and the
  corresponding highest weights $(r+1-r_i)\om_i$ will be called {\em
  pseudo fundamental weights}.
\end{defi}

Note that the pseudo fundamental weights span $P'$. By Theorem
\ref{dcrys}, the crystals of the pseudo fundamental representations of
$\Glie$ correspond to the crystals of the fundamental representations
of $^L\g$. 

Let us start with $B_2$. We have the crystal $\mathcal{M}(Y_{1,0})$
of the $5$-dimen\-sional fundamental representation of
$\U_q(B_2)$ decomposed in $\mathcal{M}^L(Y_{1,0})$ of 
the $4$-dimensional fundamental representation of
$\U_{-t}(C_2)$ and to $\mathcal{M}^L(1)$:
$$1_0\overset{1}{\rightarrow} 1_2^{-1}2_1^2\overset{2}{\rightarrow} 2_1 2_3^{-1}
\overset{2}{\rightarrow}1_22_3^{-2}\overset{1}{\rightarrow} 1_4^{-1}\text{ , }
1_0\overset{1}{\rightarrow}1_2^{-1}2_1 \overset{2}{\longrightarrow}1_22_3^{-1}\overset{1}{\rightarrow}1_4^{-2}  \text{ $\sqcup$ }\{1\}.$$
Now we have
$$
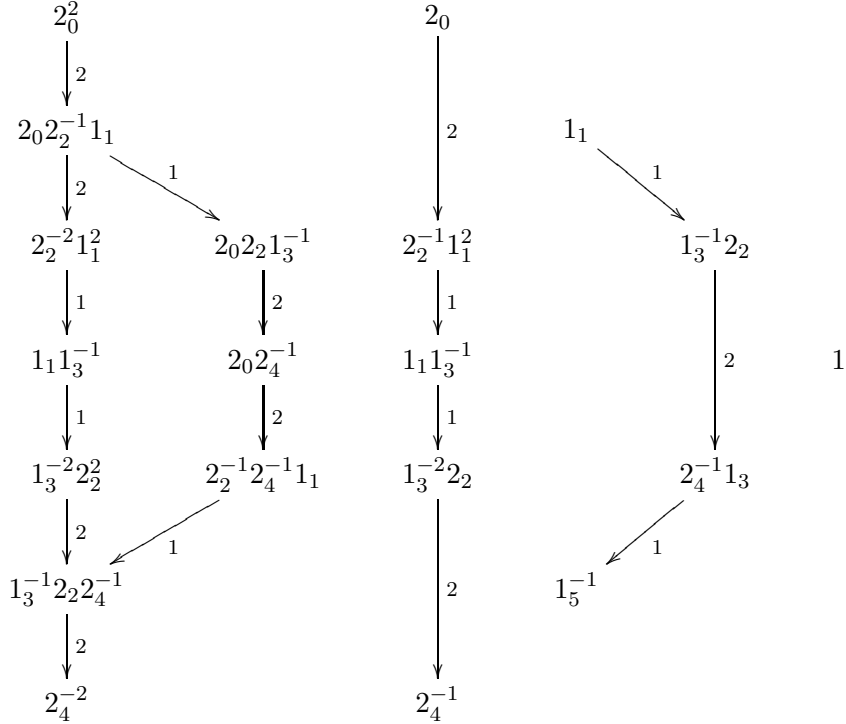
\xymatrix{2_0^2 \ar[d]^2&                          &2_0  \ar[dd]^2&        &           &                       &
\\2_0 2_2^{-1}1_1 \ar[d]^2\ar[dr]^1&                 &                       &1_1 \ar[dr]^1&                       &
\\2_2^{-2}1_1^2\ar[d]^1& 2_0 2_2 1_3^{-1}\ar[d]^2    &2_2^{-1}1_1^2 \ar[d]^1&           & 1_3^{-1}2_2\ar[dd]^2   &
\\1_1 1_3^{-1} \ar[d]^1& 2_02_4^{-1}\ar[d]^2         &1_11_3^{-1}  \ar[d]^1 &           &                       &1 
\\1_3^{-2}2_2^2\ar[d]^2&2_2^{-1}2_4^{-1}1_1\ar[dl]^1 &1_3^{-2}2_2\ar[dd]^2  &           & 2_4^{-1}1_3 \ar[dl]^1 &
\\1_3^{-1}2_22_4^{-1}\ar[d]^2&                       &                      &  1_5^{-1} &                       &  
\\2_4^{-2}                   &                       & 2_4^{-1}             &           &                       &}$$ 
The left crystal is $\mathcal{M}(Y_{2,0}^2)$ corresponding to the
$10$-dimensional representation $L(2\om_2)$ of $\U_q(B_2)$. The middle
crystal is $\mathcal{M}^L(Y_{2,0})$ corresponding to the
$5$-dimensional fundamental representation of $\U_{-t}(C_2)$. The two
right crystal contain the remaining monomials and are respectively
isomorphic to $\mathcal{M}^L(Y_{1,1})$ and $\mathcal{M}^L(1)$.

Now we suppose that $\Glie$ is of type $G_2$. $\mathcal{M}(Y_{1,0})$
has $14$ terms
$\{1_0,                        
1_2^{-1}2_1^3,                                         
2_1^22_3^{-1},                                          
\\2_12_3^{-2}1_2,                               
2_3^{-3}1_2^2, 2_12_31_4^{-1},      
1_21_4^{-1},             2_12_5^{-1},           
1_4^{-2}2_3^3,           2_3^{-1}2_5^{-1}1_2,   
1_4^{-1}2_3^22_5^{-1},                                  
2_32_5^{-2},                                          
2_5^{-3}1_4, 1_6^{-1}\}$.

The corresponding
$\tilde{\mathcal{B}}(\om_1)$ has $8$ terms 
with two connected
components described here. The first component is
$\{1_0,                        
1_2^{-1}2_1^3,                                         
2_3^{-3}1_2^2,      
1_21_4^{-1},                   
1_4^{-2}2_3^3,     
2_5^{-3}1_4, 1_6^{-1}\}$ isomorphic to $\mathcal{B}^L(\om_1)$ and
the second component is $\{2_1 2_5^{-1}\}$ isomorphic to $\mathcal{B}^L(0)$.

$\mathcal{M}(Y_{2,0}^3)$ corresponds to the $77$-dimensional
representation of $\U_q(G_2)$.  The corresponding
$\tilde{\mathcal{B}}(3\om_2)$ has $29$ terms with $4$-connected
components that we describe.

The first connected component is isomorphic to $\mathcal{B}^L(\om_2)$
($14$ terms): $\{2_0^3, 2_2^{-3}1_1^3, 1_1^21_3^{-1},
\\1_11_3^{-2}2_2^3, 1_3^{-3}2_2^6, 2_4^{-3}1_11_3, 2_2^32_4^{-3},
1_11_5^{-1}, 1_3^{-1}1_5^{-1}2_2^3, 2_4^{-6}1_3^3,
1_3^21_5^{-1}2_4^{-3}, 1_31_5^{-2}, 1_5^{-3}2_4^3, 2_6^{-3}\}$.

The second connected component is isomorphic to $\mathcal{B}^L(\om_1)$
($7$ terms): \\$\{2_02_4^{-1}1_1, 1_3^{-1}2_02_2^32_4^{-1},
2_02_4^{-4}1_3^2, 202_4^{-1}1_31_5^{-1}, 1_5^{-2}2_02_4^2,
2_2^{-1}2_6^{-2}1_1, 2_2^22_6^{-2}1_3^{-1}\}$.

The third connected component is isomorphic to $\mathcal{B}^L(\om_1)$
($7$ terms): \\$\{2_0^22_4^{-2}1_3, 2_0^22_41_5^{-1},
2_2^{-2}2_6^{-1}1_1^2, 2_22_6^{-1}1_11_3^{-1}, 2_2^42_6^{-1}1_3^{-2},
2_22_4^{-3}2_6^{-1}1_3, 1_5^{-1}2_22_6^{-1}\}$.

The fourth connected component is isomorphic to $\mathcal{B}^L(0)$
($1$ term): $\{2_02_22_4^{-1}2_6^{-1}\}$.

Now let us look at the remaining example of Section \ref{charstate} for $\Glie$ is of type $B_2$
$$\xymatrix{1_0^2 \ar[d]^1&                                     &1_0^2 \ar[d]^1&                         &&
\\1_0 1_2^{-1}2_1^2\ar[d]^1\ar[dr]^2&                           &1_0 1_2^{-1}2_1\ar[d]^1\ar[ddr]^2&       &&
\\1_2^{-2}2_1^4\ar[d]^2&1_0 2_1 2_3^{-1}\ar[dl]^1\ar[d]^2       &1_2^{-2}2_1\ar[dd]^2&                &&1_0\ar[dl]^1
\\1_2^{-1}2_1^32_3^{-1} \ar[d]^2& 1_01_22_3^{-2}\ar[d]^1        && 1_01_22_3^{-1}\ar[d]^1       &1_2^{-1}2_1\ar[dd]^2&
\\2_1^22_3^{-2}\ar[d]^2&1_01_4^{-1}\ar[d]^1                     &2_12_3^{-1}\ar[dd]^2&1_01_4^{-1}\ar[d]^1 &&
\\2_12_3^{-3}1_2\ar[d]^2\ar[dr]^1&1_2^{-1}1_4^{-1}2_1^2\ar[d]^2 &&1_2^{-1}1_4^{-1}2_1\ar[ddl]^2         &2_3^{-1}1_2\ar[dr]^1
\\2_3^{-4}1_2^2\ar[d]^1&1_4^{-1}2_12_3^{-1}\ar[dl]^2            &2_3^{-2}1_2^2\ar[d]^1&              &&1_4^{-1} 
\\2_3^{-2}1_21_4^{-1}\ar[d]^1&                                  &2_3^{-1}1_21_4^{-1}\ar[d]^1&            &&
\\1_4^{-2}&                                                     &1_4^{-2}&                               &&}$$ 
The left crystal is $\mathcal{M}(Y_{1,0}^2)$ corresponding to the
$14$-dimensional representation $L(2\om_1)$ of $\U_q(B_2)$. The middle
crystal is $\mathcal{M}^L(Y_{1,0}^2)$ corresponding to the
$10$-dimensional representation $L(2\om_1)$ of $\U_{-t^2}(C_2)$. The
right crystal contains the remaining monomials and is isomorphic to
$\mathcal{M}^L(Y_{1,0})$.

\subsection{Proof of Theorem \ref{dcrys}}

We consider operators $f_i^L, e_i^L$ on $\mathcal{M}$
as defined in formula \eqref{tilde f e}. Let
\begin{align*}
\mathcal{M}' &= \{m'\in\mathcal{M}|\forall i \in I, l\in\ZZ; u_{i,l} \in
(r+1-r_i)\ZZ \} \\
&= \{m'\in\mathcal{M}|\forall i \in I, l\in\ZZ; u_{i,l}\in
r\ZZ\text{ if }r_i = 1\}.
\end{align*}
As $\mathrm{wt}(\mathcal{M}')\subset P'$, we can define $\mathrm{wt}^L
= \Pi\circ\mathrm{w}$ on $\mathcal{M}'$.

\begin{lem}\label{aidecris} Let $i\in I$ such that $r_i = 1$. Let
  $m\in\mathcal{M}'$ such that $\epsilon_i(m) > 0$ (resp. $\phi_i(m) >
  0$). Then for $1\leq q\leq r - 1$ we have

(1) $\epsilon_i(e_i^q(m)) > 0$
    (resp. $\phi_i(f_i^q(m)) > 0$),

(2) $p_i(e_i^q(m)) = p_i(m)$ (resp. $q_i(f_i^q(m)) =
    q_i(m)$),

(3) $e_i^r(m) = m A_{i,p_i(m) - 1}^r$ (resp. $f_i^r(m)
    = m A_{i,q_i(m) + 1}^{-r}$).
\end{lem}

\demo We prove the assertions for $\phi_i(m) > 0$ (the assertions for
$\epsilon_i(m) > 0$ are proved in the same way).

As $m\in\mathcal{M}'$, we have $\phi_i(m) \in r\ZZ$, and so $\phi_i(m)
\geq r$. So $\phi_i(f_i^q(m)) = \phi_i(m) - q \geq 1$ and the
statement (1) is proved.

We have $f_i(m) = m A_{i,q_i(m) + 1}^{-1}$. We have
$\phi_{i,q_i(m)} (f_i(m)) = \phi_i(m) - 1$. For $l\geq q_i(m)
+ 2$, we have $\phi_{i,l}(f_i(m)) = \phi_{i,l}(m) - 2\leq
\phi_i(m) - 2$. For $l < q_i(m)$, we have $\phi_{i,l}(m)\in r\ZZ$, so
$\phi_{i,l}(m)\leq \phi_i(m) - r$ and $\phi_{i,l}(f_i(m)) =
\phi_{i,l}(m) \leq \phi_i(m) - r$. So $q_i(f_i(m)) = q_i(m)$
and we have proved the point (2) for $q = 1$. If $r = 3$ we also have
to prove the statement for $q = 2$. We have $f_i^2(m) = m
A_{i,q_i(m) + 1}^{-2}$. We have $\phi_{i,q_i(m)} (f_i(m)^2) =
\phi_i(m) - 2$. For $l\geq q_i(m) + 2$, we have
$\phi_{i,l}(f_i^2(m)) = \phi_{i,l}(m) - 4\leq \phi_i(m) -
4$. For $l < q_i(m)$, we have $\phi_{i,l}(m)\in r\ZZ$, so
$\phi_{i,l}(m)\leq \phi_i(m) - r$ and $\phi_{i,l}(f_i^2(m)) =
\phi_{i,l}(m) < \phi_i(m) - 2$. So $q_i(f_i(m)) = q_i(m)$ and
we have proved the point (2) for $q = 2$.

The last assertion (3) is a direct consequence of the first two
assertions.
\qed

Let $\Pi : \mathcal{M}'\rightarrow \mathcal{M}$ be the map defined by 
$$\Pi(m) = \prod_{i\in I}Y_{i,l}^{u_{i,l}(m)(1 + r - r_i)^{-1}}.$$
Let $\mathcal{M}^L$ be the monomial crystal for $^L\g$. Viewed as a
set, $\mathcal{M}^L$ is equal to $\mathcal{M}$ and so we can view the
map $\Pi$ as $\Pi : \mathcal{M}'\rightarrow \mathcal{M}^L$.

\begin{thm}\label{ld} $\mathcal{M}'\sqcup \{0\}$ is stable for the
  operators $f_i^L$, $e_i^L$ which define a structure
  of $^L\g$-crystal on $\mathcal{M}'$. The map $\Pi :
  \mathcal{M}'\rightarrow \mathcal{M}^L$ is an isomorphism of
  $^L\g$-crystals.
\end{thm}

\demo The stability for $e_i,f_i$ when $r_i = r$ is
clear as the $A_{i,l}^{\pm}\in\mathcal{M}'$. When $r_i = 1$ it is a
consequence of Lemma \ref{aidecris} as the $A_{i,l}^{\pm
  r}\in\mathcal{M}'$.

To prove that we have a crystal isomorphism, first note that the
compatibility of the map with $\epsilon_i, \phi_i$ is clear. Then for
the compatibility with the operators $e_i, f_i$, it is
clear if $r_i = r$ and if $r_i = 1$ it follows from Lemma
\ref{aidecris}.
\qed

Theorem \ref{dcrys} is a direct consequence of Theorem \ref{ld}. Thus, 
Theorem \ref{dcrys} is now proved. As discussed above Theorem \ref{incchar}
is also now proved.

\begin{rem} The proof given above also implies that Theorems \ref{dcrys} and \ref{incchar} hold for any symmetrizable Kac--Moody algebra  such that $r\leq 3$ (see \cite{hn} for their monomial crystal). 
Here the $r_i$ are defined as the set of relatively prime integers such that $r_i C_{i,j} = r_j C_{j,i}$, and $r$ is the maximal
number among the $r_i$.
\end{rem}

\section{Interpolating quantum groups}    \label{iqg}

In the previous section we have described a duality between characters
and crystal bases of finite-dimensional representations of $U_q(\g)$
and $U_q({}^L\g)$. We would like to explain this duality in the
following way: there exists a two-parameter deformation of both of
these quantum groups, which we call the ``interpolating quantum
group''. Moreover, the dual finite-dimensional representations
$U_q(\g)$ and $U_q({}^L\g)$ appear as the result of specialization (of
the first and the second parameter, respectively) of a representation
of this interpolating quantum group.

In this section we define the interpolating quantum group and in the
following two sections we construct their representations which
exhibit the desired duality property.

Let again $\g$ be a finite-dimensional simple Lie algebra and
$U_q(\g)$ the corresponding quantum group. We denote by $\U_q(\Glie)$
the algebra with the same generators and relations except for the
Serre relations. Note that $U_q(\g)$ and $\U_q(\Glie)$ have the same
categories of finite-dimensional representations.

The interpolating quantum group $\U_{q,t}(\g)$ is an associative
algebra depending on two parameters, $q$ and $t$. (Note that this
algebra is different from the two-parameter quantum groups considered
in \cite{bw, r}.) We will then establish the following Langlands
duality property of these algebras: the specialization with respect to
one parameter, $t=1$, gives the quantum group $\U_q(\g)$, and the
specialization with respect to the other parameter, $q=\ep$, where
$\ep = 1$ for simply-laced $\g$ and $\exp(\pi i/r)$ for non-simply
laced ones, gives the Langlands dual quantum group $\U_{-t}({}^L\g)$.

\subsection{Interpolating simply-laced quantum groups ($r = 1$)}

Let $\Glie$ be a simply-laced simple Lie algebra, that is, $r=1$. In
this situation the definition of the interpolating quantum group is
essentially equivalent to the usual definition of quantum group. In
what follows by an ``algebra'' we will always mean an associative
unital algebra over $\C$.

\begin{defi}
$\U_{q,t}(\Glie)$ is the algebra with the generators $X_i^{\pm}$,
$K_i^{\pm 1}$, $\wt{K}_i^{\pm 1}$ and relations
$$K_i X_j^{\pm} = q^{\pm C_{i,j}}X_j^{\pm}K_i\text{ , }\wt{K_i}
X_j^{\pm} = t^{\pm C_{i,j}}X_j^{\pm}\wt{K}_i,$$
$$[X_i,X_j^-] = \delta_{i,j} \frac{K_i\wt{K}_i -
  (K_i\wt{K}_i)^{-1}}{qt - (qt)^{-1}}.$$
\end{defi}

Note that
$$
\U_{q,t}(\Glie) \supset \, \langle (K_i\tilde{K_i})^{\pm 1},X_i^\pm
\rangle \, \simeq \U_{qt}(\Glie),
$$
and that we have the following interpolating property:
$$
\U_{q,1}(\Glie)/(\tilde{K}_i = 1)\simeq \U_q(\Glie) \qquad \text{
and } \qquad \U_{1,t}(\Glie)/(K_i = 1) = \U_t(\Glie) =
\U_t({}^L\g).
$$
As a special case, we have the elementary interpolating quantum group
$\U_{q,t}(A_1)$. The elementary rank one subalgebras of
$\U_{q,t}(\Glie)$ corresponding to simple roots are all isomorphic to
$\U_{q,t}(A_1)$ if $\Glie$ is simply-laced. This is analogous to
the properties of standard quantum groups. We will see in the
following that for non-simply laced $\Glie$ we will have to consider
other elementary (rank $1$) interpolating quantum groups
corresponding to $B_1$, $C_1 = {}^L B_1$, $G_1$, and $^L G_1$.

\subsection{Elementary interpolating quantum groups for $r = 2$}

For $r = 2$ we have $\ep = \exp(\pi i/2) = i$. We will define two
elementary interpolating quantum groups $\U_{q,t}(C_1)$ and
$\U_{q,t}(B_1)$. The definition of the first one is simple.

\begin{defi} $\U_{q,t}(C_1)$ is the algebra with generators $X^\pm$,
  $K^{\pm 1}$, $\tilde{K}^{\pm 1}$ and relations
$$KX^{\pm} = q^{\pm 4} X^{\pm}K\text{ , }\tilde{K}X^{\pm} = t^{\pm 2}
X^{\pm}\tilde{K},$$
$$[X^+,X^-] = \frac{K\tilde{K} - (K\tilde{K})^{-1}}{q^2t - q^{-2}t^{-1}}.$$
\end{defi}

Note that $$\U_{q,t}(C_1) \supset \, \langle (K\tilde{K})^{\pm
  1},X^\pm \rangle \, \simeq
\U_{q^2t}(sl_2),$$ and that we have the interpolating property
$$\U_{q,1}(\Glie)/(\tilde{K} = 1)\simeq \U_{q^2}(sl_2) =
\U_q(C_1)\text{ and }\U_{\epsilon,t}(\Glie)/(K = 1) \simeq \U_{-t}(sl_2) =
\U_{-t}({}^L C_1),$$ as $\U_{-t}(sl_2) = \U_{-t}(B_1)$.

\begin{defi} $\U_{q,t}(B_1)$ is the algebra with generators $X^{\pm
  }$, $K^{\pm 1}$, $\tilde{K}^{\pm 1}$, $\eta$, central elements
  $C$, $\tilde{C}$ and relations
$$\CC[K^{\pm 1},\tilde{K}^{\pm 1},\eta] \qquad \text{is
    commutative},$$
$$K X^\pm = q^{\pm 2}X^\pm K\text{ , }\tilde{K}X^\pm = t^{\pm 1}X^\pm
  \tilde{K} \text{ , } \eta X^\pm = X^\pm (\eta \pm 1),$$
\begin{equation}\label{hun} X^\pm X^\mp = \frac{q^C (t^{\tilde{c}}
\tilde{K}^{\pm 1})^P + q^{-C}(t^{\tilde{c}}\tilde{K}^{\pm 1})^{-P} -
q^{\mp 1}t^{\pm \tilde{c}}\tilde{K} K - q^{\pm 1} t^{\mp \tilde{c}}
(\tilde{K}K)^{-1}}{(q-q^{-1})(qt - (qt)^{-1})},\end{equation}
where $P = (-1)^\eta$ and $\tilde{c} = P \tilde{C} - 1/2$.
\end{defi}

Note that we have $t^{\tilde{c}}X^\pm = X^\pm t^{-\tilde{c}-1}$, $P^2
= 1$ and $P$ commutes with $E^2$ and $F^2$. We also have the following:
$$q^C (t^{\tilde{c}}\tilde{K}^{\pm 1})^P +
q^{-C}(t^{\tilde{c}}\tilde{K}^{\pm 1})^{-P} = q^{PC}
t^{\tilde{c}}\tilde{K}^{\pm 1} + q^{-PC}t^{-\tilde{c}}\tilde{K}^{\mp
1}.$$
The elements
$$Cas(q) = q^C+q^{-C}\text{ and }Cas(t) = t^{\tilde{c}+
  1/2}+t^{-\tilde{c} - 1/2}$$ are central. The element $Cas(q)$ will
  correspond to the Casimir element for the specialization $t=1$. For
  the other specialization, $q=\ep$, the Casimir element will not
  be exactly $Cas(t)$, but $t^{2(1 + \tilde{c})} + t^{-2(1+\tilde{c})}$,
  which is not central in the whole algebra, but commutes with $(X^\pm
  )^2$.

\begin{lem} The algebra $\U_{q,t}(B_1)$ is well-defined.\end{lem}

\demo The only point to be checked is the associativity condition
$(X^\pm X^\mp)X^\pm = X^\pm (X^\mp X^\pm)$. It is satisfied as we have
$$(q^C(t^{\tilde{c}}\tilde{K}^{\pm 1})^P + q^{-C}(t^{\tilde{c}}
\tilde{K}^{\pm 1})^{-P} - q^{\mp 1}t^{\pm \tilde{c}}\tilde{K} K -
q^{\pm 1} t^{\mp \tilde{c}} (\tilde{K}K)^{-1})X^\pm $$
$$=X^\pm(q^C(t^{-\tilde{c}}\tilde{K}^{\pm 1})^{-P} +
q^{-C}(t^{-\tilde{c}}\tilde{K}^{\pm 1})^P - q^{\pm 1}t^{\mp
\tilde{c}}\tilde{K} K - q^{\mp 1} t^{\pm \tilde{c}}
(\tilde{K}K)^{-1}).$$ \qed

Let us look at the specializations of $\U_{q,t}(B_1)$ at $t=1$ and
$q=\ep=i$. Let
$$\mathcal{X}^\pm = \mp (X^\pm)^2/(t - t^{-1})\text{ , }
\mathcal{K} = \tilde{K}^2.$$ 

\begin{prop}    \label{inter}
The subalgebra of $\U_{q,1}(B_1)/(\tilde{K} = 1)$
generated by $X^\pm$, $K^{\pm 1}$ is isomorphic to $\U_q(sl_2) =
\U_q(B_1)$.

The subalgebra of $\U_{\ep,t}(B_1)/(K^2 = 1, Kq^{PC} = \ep)$ generated by
$\mathcal{X}^\pm$, $\mathcal{K}^{\pm 1}$ is isomorphic to
$\U_{t^2}(sl_2) = \U_{-t}(C_1) = \U_{-t}({}^L B_1)$.
\end{prop}

\demo First, let us consider the specialization $\U_{q,1}(B_1)$ at $t
= 1$. Then the element $\tilde{K}$ becomes central and we can 
specialize $\tilde{K} = 1$. We have the
relations $KX^\pm = q^{\pm 2}X^\pm K$ and
$$(q-q^{-1})^2 X^\pm X^\mp + q^{\mp 1}K + q^{\pm 1} K^{-1} = Cas(q).$$ The equality implies the standard relation
$$[X^+,X^-] = \frac{K - K^{-1}}{q - q^{-1}}.$$ $Cas(q)$ is central and
corresponds to the central Casimir element in $\U_q(sl_2)$. So we have
an isomorphism.

Now let us consider the specialization of $\U_{q,t}(B_1)$ at $q =
\ep$. Then $K^2$ becomes central. Let us consider the algebra
$\U_{\ep,t}(B_1)/(K^2 = 1)$. We have the relations:
$$
KX^\pm = -X^\pm K\text{ , }KX^\pm = - X^\pm K\text{ ,
}\tilde{K}X^\pm = t^{\pm 1}X^\pm \tilde{K},$$
$$X^{\pm} X^{\mp} =\frac{(q^{PC} + \ep K)(t^{\tilde{c}}\tilde{K}^{\pm 1}
- \ep q^{-PC}K t^{-\tilde{c}} \tilde{K}^{\mp 1})}{-2(t + t^{-1})}.$$
Since $(q^{PC} + \ep K)X^\pm = X^\pm(q^{-PC} - \ep K) = X^\pm(-\ep
Kq^{-PC})(q^{PC}+ \ep K)$,
we find that $4(t+t^{-1})^2(X^\pm)^2(X^\mp)^2$ is equal to
$$
-\ep K q^{-PC}(q^{PC}K + \ep)^2(t^{-\tilde{c}-2}\tilde{K}^{\pm 1} +
\ep q^{PC}K t^{\tilde{c}+2} \tilde{K}^{\mp 1})
(t^{\tilde{c}}\tilde{K}^{\pm 1} - \ep q^{-PC}K t^{-\tilde{c}}
\tilde{K}^{\mp 1}).
$$
So it is natural to specialize at $K q^{PC} = \ep$. We obtain that
$$(X^\pm)^2(X^\mp)^2 
= \frac{t^{-2} \tilde{K}^{\pm 2} + t^2 \tilde{K}^{\mp 2} -
  t^{2(1+\tilde{c})} - t^{-2(1+\tilde{c})}}{(t+t^{-1})^2}.$$
The above relations can be rewritten as
$$ (t^2 - t^{-2})^2\mathcal{X}^{\pm}\mathcal{X}^{\mp} + t^{\mp
2}\mathcal{K}+t^{\pm 2}\mathcal{K}^{-1} = t^{2(1+\tilde{c})} + t^{-2(1+\tilde{c})}.$$
The element $t^{2(1+\tilde{c})} +
t^{-2(1+\tilde{c})}$ commutes with $\mathcal{X}^{\pm}$,
$\mathcal{K}^{\pm 1}$ and corresponds to the Casimir element (see the
above discussion).
We get the equality
$$[\mathcal{X}^+,\mathcal{X}^-] = \frac{\mathcal{K} -
\mathcal{K}^{-1}}{t^2 - t^{-2}}.$$
\qed

\subsection{Interpolating quantum group for $r = 2$} Let $\Glie$ be a
simple Lie algebra such that $r = 2$, that is, $\Glie$ is of type
$B_n$, $C_n$ or $F_4$.

\begin{defi} $\U_{q,t}(\Glie)$ is the algebra with generators
  $X_i^\pm$, $K_i^{\pm 1}$, $\tilde{K}_i^{\pm 1}$, $\eta_j$,
  $C_j$, $\tilde{C}_j$ ($1\leq i,j\leq n$, $r_j = 1$) and relations
$$\CC[K_i^{\pm},\tilde{K}_i^{\pm 1},\eta_j, C_j,
\tilde{C}_j]_{1\leq i,j\leq n,r_j = 1} \qquad \text{is commutative},$$
$$\mathcal{U}_i = \langle X_i^\pm,K_i^{\pm 1},\tilde{K}_i^{\pm 1}
  \rangle \simeq \U_{q,t}(C_1)\text{ if $r_i = 2$,}$$
$$\mathcal{U}_i = \langle X_i^\pm, K_i^{\pm 1}, \tilde{K}_i^{\pm 1},
\eta, C_i, \tilde{C}_i \rangle \simeq \U_{q,t}(B_1)\text{ if $r_i
= 1$,}$$
$$K_i X_j^\pm = q^{\pm r_i C_{i,j}} X_j^\pm K_i\text{ , }\tilde{K}_i
X_j^\pm = t^{\pm r_iC_{i,j}/2}\tilde{K}_j,$$
$$[X_i^+ , X_j^-] = [(-1)^{\eta_i} , X_j^\pm] = 0\text{ for $i\neq
  j$.}$$
\end{defi}

Let us consider the elements
$$\mathcal{X}_i^\pm =\mp (X_i^\pm)^2/(t - t^{-1})\text{ ,
}\mathcal{K}_i = \tilde{K}_i^2\text{ for $r_i = 1$,}$$
$$\mathcal{X}_i^\pm =X_i^\pm\text{ , }\mathcal{K}_i =
\tilde{K}_i\text{ for $r_i = 2$.}$$

The specialization at $q = \ep = i = \sqrt{-1}$ should not be confused
in the following with the index $i\in I$. Proposition \ref{inter}
implies

\begin{prop}    \label{inter1}
The subalgebra of $\U_{q,1}(\Glie)/(\tilde{K}_i = 1)$ generated by the
  $X_i^\pm$, $K_i^{\pm 1}$ is isomorphic to $\U_q(\Glie)$.

The subalgebra of $\U_{\ep,t}(\Glie)/(K_i^2 = 1, K_iq^{P_iC_i} = \ep)$
generated by the $\mathcal{X}_i^\pm$, $\mathcal{K}_i^{\pm 1}$ is
isomorphic to $\U_{-t}({}^L\g)$.
\end{prop}

In the proposition, by convention, $P_iC_i = 1$ if $r_i = 2$, that is,
the relation $K_iq^{P_iC_i} = \ep$ means $K_i = 1$.

According to the above proposition, $\U_{q,t}(\Glie)$ interpolates
between $\U_q(\Glie)$ and $\U_{-t}({}^L\g)$ the quantum groups without
the Serre relations. Is it possible to have an algebra that
interpolates between the quantum groups $U_q(\Glie)$ and
$U_{-t}({}^L\g)$ with the Serre relations? In other words, can one
construct a two-parameter deformation of the Serre relations of
$U_q(\Glie)$ and $U_{-t}({}^L\g)$? In this paper we are only
interested in finite-dimensional representations. Therefore this
question is not important, because finite-dimensional representations
of $\U_q(\Glie)$ are the same as those of $\U_q(\Glie)$ (and similarly
for $\U_{-t}({}^L\g)$ and $U_{-t}({}^L\g)$). But for other
representations this question becomes important. The examples given
below indicate that in the framework of $\U_{q,t}(\Glie)$ the answer
is negative.

In fact, in Example 3 of Section \ref{int reps} we will construct a
finite-dimensional representation $V$ of $\U_{q,t}(B_2)$ which
interpolates between representations of $\U_q(B_2)$ and $\U_{-t}(C_2)$
(and hence of $U_q(B_2)$ and $U_{-t}(C_2)$), but for different vectors
in this representation different $t$-deformations of the Serre
relations of $U_q(B_2)$ will be satisfied. Imposing either of them (or
another $t$-deformation) on the algebra would lead to additional
relations that are not satisfied in $V$. Therefore $V$ is not a module
over this algebra. Hence it appears impossible to incorporate a
two-parameter deformation of the Serre relations into
$\U_{q,t}(\Glie)$ in such a way that Proposition \ref{inter1} would
hold for the quotient, with $\U_q(\Glie)$ and $\U_{-t}({}^L\g)$
replaced by $U_q(\Glie)$ and $U_{-t}({}^L\g)$.

To illustrate this point further, consider the following example of a
candidate for a $t$-deformation of the Serre relations for $\g=B_2$
(note that we do not use it in this paper):
\begin{equation}    \label{faux serre}
X^+_2 {X^+_1}^2 - (q^2 t + q^{-2} t^{-1}) X^+_1X^+_2X^+_1+ {X^+_1}^2
X^+_2 = 0,
\end{equation}
$$X^+_1{X^+_2}^3 - t(q^2 + 1 + q^{-2}) X^+_2X^+_1{X^+_2}^2 +
t^{-2}(q^2 + 1 + q^{-2}) {X^+_2}^2X^+_1X^+_2 - t^{-1}{X^+_2}^3X^+_1 =
0.$$ At $t=1$ we recover the Serre relation of $\U_q(B_2)$. Let us
consider the specializations $S$, $S'$ of these relations at $q = \ep$.
By computing $S'X^+_2 - t X^+_2S'$ we obtain
$$X^+_1({X^+_2}^2)^2  - (t^2 + t^{-2}) ({X^+_2}^2)X^+_1({X^+_2}^2) +
({X^+_2}^2)^2X^+_1 = 0,$$
which is one of the Serre relation of $\U_{-t}(C_2)$. By computing 
$$X^+_2 S X^+_1 - (t^2 + 1 + t^{-2}) X^+_1X^+_2S + (t^2 + 1 + t^{-2})
S X^+_2X^+_1$$ $$- X^+_1 S X^+_2 - (t + t^{-1}) X^+_2X^+_1 S + (t +
t^{-1}) S X^+_1X^+_2,$$ we obtain
$$0 = ({X^+_2}^2){X^+_1}^3 -(t^2 + 1 + t^{-2})
X^+_1({X^+_2}^2){X^+_1}^2 + (t^2 + 1 + t^{-2})
{X^+_1}^2({X^+_2}^2)X^+_1 - {X^+_1}^3({X^+_2}^2),$$ which is another
Serre relation for $\U_{-t}(C_2)$ (both relations should be written in
terms of $X_1^+$ and $\mathcal{X}_2^+ = (X_2^+)^2$).

But if we compute the bracket of the second Serre relation with
$X^-_1$, we obtain
$$K_1{X^+_2}^3q^2(1 + q^2)(t^{-4} - 1) = (1-t^2)(1-q^{-2})^2(1 +
q^{-2}) K_1^{-1}{X^+_2}^3.$$
Then we following identity which does not hold in either $U_q(B_2)$ or
$\U_{-t}(C_2)$:
$$K_1{X^+_2}^3 q^4 (1 + t^2)=t^4(1 - q^{-2})^2 K_1^{-1} {X^+_2}^3.$$
Hence if we include the relations \eqref{faux serre}, we obtain an
algebra that does not have the desired interpolation property.

\subsection{Interpolating quantum groups for $r = 3$} For $r = 3$ we
define two elementary interpolating quantum groups $\U_{q,t}(G_1)$ and
$\U_{q,t}({}^LG_1)$. We have $\epsilon = e^{2\pi i/6}$.

\begin{defi} $\U_{q,t}({}^LG_1)$ is the algebra with generators $X^\pm$,
  $K^{\pm 1}$, $\tilde{K}^{\pm 1}$ and relations
$$KX^{\pm} = q^{\pm 6} X^{\pm}K\text{ , }\tilde{K}X^{\pm} = t^{\pm 2}
X^{\pm}\tilde{K},$$
$$[X^+,X^-] = \frac{K\tilde{K} - (K\tilde{K})^{-1}}{q^3t -
  q^{-3}t^{-1}}.$$
\end{defi}

Note that
$$ \U_{q,t}({}^LG_1) \supset \, \langle (K\tilde{K})^{\pm
1},X^\pm\rangle \, \simeq \U_{q^3t}(sl_2)$$ and that we have the
following interpolating property:
$$\U_{q,1}(\Glie)/(\tilde{K} = 1)\simeq \U_{q^3}(sl_2) =
\U_q({}^LG_1)\text{ and }\U_{\epsilon,t}(\Glie)/(K = 1) \simeq \U_{-t}(sl_2) =
\U_{-t}(G_1).$$

Let us define the elementary interpolating quantum group
$\U_{q,t}(G_1)$. First we need the following polynomial map $F(X) =
X(X - \epsilon^2)(1 - \epsilon^2)^{-1}$ which satisfies $F(1) = 1$,
$F(\epsilon^2) = 0$, $F(\epsilon^4) = -1$.

\begin{defi} We define the algebra $\U_{q,t}(G_1)$ as the algebra with
  generators $X^{\pm}$, $K^{\pm 1}$, $\tilde{K}^{\pm 1}$, $\eta$,
  central elements $C, \tilde{C}$ and relations
$$\CC[K^{\pm 1},\tilde{K}^{\pm 1},\eta] \quad \text{is
    commutative},$$
$$K X^\pm = q^{\pm 2}X^\pm K\text{ , } \tilde{K}X^\pm = t^{\pm 1}X^\pm
\tilde{K}
\text{ , }\eta X^\pm = X^\pm (\eta \pm 1),$$
\begin{multline}\label{hung} X^\pm X^\mp = \\
  q^C(t^{\tilde{c}_\pm}\tilde{K})^{P_\pm} + q^{-C}
  (t^{\tilde{c}_\pm}\tilde{K})^{-P_\pm}
- q^{\mp 1}K (t^{\tilde{c}_\pm}\tilde{K})^{P_\pm^2} - q^{\pm
  1}K^{-1}(t^{\tilde{c}_\pm}\tilde{K})^{-P_\pm^2},\end{multline}
where $P_\pm = F(\epsilon^{2(-\eta + 1 \mp 1)})$ and $\tilde{c}_\pm = P_\pm
  \tilde{C} \mp 1/2$.
\end{defi}

\begin{lem} The algebra $\U_{q,t}(G_1)$ is well-defined.\end{lem}

\demo The only point to be checked is the associativity condition
$(X^\pm X^\mp)X^\pm = X^\pm (X^\mp X^\pm)$, which is verified as
follows: $$(q^C(t^{\tilde{c}_\pm}\tilde{K})^{P_\pm} + q^{-C}
(t^{\tilde{c}_\pm}\tilde{K})^{-P_\pm} - q^{\mp 1}K
(t^{\tilde{c}_\pm}\tilde{K})^{P_\pm^2} - q^{\pm
1}K^{-1}(t^{\tilde{c}_\pm}\tilde{K})^{-P_\pm^2}) X^\pm$$ $$=X^\pm
(q^C(t^{\tilde{c}_\mp \mp 1}\tilde{K}t^{\pm 1})^{P_\mp} + q^{-C}
(t^{\tilde{c}_\mp \mp 1}\tilde{K}t^{\pm 1})^{-P_\mp}$$ $$- q^{\pm 1}K
(t^{\tilde{c}_\mp\mp 1}\tilde{K}t^{\pm 1})^{P_\mp^2} - q^{\mp
1}K^{-1}(t^{\tilde{c}_\mp\mp 1}\tilde{K}t^{\pm 1})^{-P_\mp^2})
$$ $$=X^\pm (q^C(t^{\tilde{c}_\mp}\tilde{K})^{P_\mp} + q^{-C}
(t^{\tilde{c}_\mp}\tilde{K})^{-P_\mp}- q^{\pm 1}K
(t^{\tilde{c}_\mp}\tilde{K})^{P_\mp^2} - q^{\mp
1}K^{-1}(t^{\tilde{c}_\mp}\tilde{K})^{-P_\mp^2}).$$ \qed

Let us set
$$\mathcal{X}^+ = \frac{(X^+)^3}{(1 - \epsilon^4)^2(t^3 - t^{-3})}\text{ ,
}\mathcal{X}^- = \frac{(X^-)^3(-1)^{m+1}}{(1 + \epsilon^4 +
  2\epsilon^5)(t^{-3} - t^3)}\text{ , }\mathcal{K} = \tilde{K}^2.$$

\begin{prop}\label{speg} The subalgebra of $\U_{q,1}(G_1)/(\tilde{K} =
  1)$ generated by $X^\pm/(q-q^{-1})$, $K^{\pm 1}$ is isomorphic to
  $\U_q(sl_2) = \U_q(G_1)$.

For any $m \in \Z/2\Z$, the quotient by $\epsilon^{2\eta} = 1$, $K =
(-1)^m$, $q^C = (-1)^m\epsilon$ of the subalgebra of
$\U_{\epsilon,t}(G_1)$ generated by $\mathcal{X}^\pm$, $\mathcal{K}$
is isomorphic to $\U_{-t^3}(sl_2) = \U_{-t}({}^LG_1)$.
\end{prop}

\demo The first point is proved as for $\U_{q,1}(B_1)$. 
Now let us consider the specialization of $\U_{q,t}(G_1)$ at $q
= \epsilon$. Then $K^3$ becomes central. 
Note that we have, $P_\pm X^\mp= X^\mp P_0$, where $P_0
= F(\epsilon^{2 - 2\eta})$.  We also have $P_0X^\pm = X^\pm P_\pm$ and
$$\tilde{c}_\pm X^\mp = X^\mp (P_0 \tilde{C} \mp 1/2)\text{ , }(P_0
\tilde{C} \mp 1/2)X^\mp = X^\mp (\tilde{c}_\mp \mp 1).$$ So we can
compute $(X^\pm)^3(X^\mp)^3$, and we obtain
$$(q^C(t^{\tilde{c}_\pm}\tilde{K})^{P_\pm} + q^{-C}
(t^{\tilde{c}_\pm}\tilde{K})^{-P_\pm} - q^{\mp 1}K
(t^{\tilde{c}_\pm}\tilde{K})^{P_\pm^2} - q^{\pm
1}K^{-1}(t^{\tilde{c}_\pm}\tilde{K})^{-P_\pm^2})$$
$$\times(q^C(t^{P_0\tilde{C} \mp 3/2}\tilde{K})^{P_0} + q^{-C}
(t^{P_0\tilde{C} \mp 3/2}\tilde{K})^{-P_0}$$ $$- q^{\mp 3}K
(t^{P_0\tilde{C} \mp 3/2}\tilde{K})^{P_0^2} - q^{\pm
3}K^{-1}(t^{P_0\tilde{C} \mp 3/2}\tilde{K})^{-P_0^2})$$
$$\times(q^C(t^{\tilde{c}_\mp \mp 3}\tilde{K})^{P_\mp} + q^{-C}
(t^{\tilde{c}_\mp \mp 3}\tilde{K})^{-P_\mp} - q^{\mp 5}K
(t^{\tilde{c}_\mp \mp 3}\tilde{K})^{P_\mp^2} - q^{\pm
5}K^{-1}(t^{\tilde{c}_\mp \mp 3}\tilde{K})^{-P_\mp^2}).$$

Note that $\epsilon^{2\eta}$ and $K$ commute with $(X^\pm)^3$. So we may
take the quotient by the relations $K = (-1)^m$, $\epsilon^{2\eta} = 1$ and
$q^C = (-1)^m\epsilon$ where $m\in\ZZ/2\ZZ$. In particular, we have $P_0 =
0$, $P_\pm = \pm 1$, $\tilde{c}_\pm = \pm\tilde{C}\mp 1/2$. For
$(X^\pm)^3(X^\mp)^3$ we obtain
$$(-1)^{m+1}(1 + \epsilon^4 + 2 \epsilon^5)(t^{\pm
  \tilde{c}_\pm}\tilde{K}^{\pm 1} + \epsilon^4
  t^{\mp\tilde{c}_\pm}\tilde{K}^{\mp 1} - \epsilon^{2 \pm 2}
  t^{\mp\tilde{c}_\pm}\tilde{K}^{\pm 1} - \epsilon^{ 2\mp 2}t^{\pm
  \tilde{c}_\pm}\tilde{K}^{\mp 1})
$$ $$\times(t^{\mp\tilde{c}_\mp + 3}\tilde{K}^{\mp 1} + \epsilon^4
  t^{\pm \tilde{c}_\mp - 3}\tilde{K}^{\pm 1} - \epsilon^{2 \mp 2}
  t^{\tilde{c}_\mp \mp 3}\tilde{K} - \epsilon^{2 \pm 2} t^{\mp 3 -
  \tilde{c}_\mp }\tilde{K}^{-1})
$$ $$= (-1)^{m+1}(1 + \epsilon^4 + 2 \epsilon^5)(1 -
\epsilon^4)^2(t^{\mp\tilde{c}_\mp}\tilde{K}^{\pm 1} -
t^{\mp\tilde{c}_\pm}\tilde{K}^{\mp 1}) (t^{\mp\tilde{c}_\mp +
3}\tilde{K}^{\mp 1} - t^{\pm \tilde{c}_\mp - 3}\tilde{K}^{\pm 1})
$$ $$= (-1)^{m+1}(1 + \epsilon^4 + 2 \epsilon^5)(1 -
\epsilon^4)^2(-t^{-3}\tilde{K}^{\pm 2} - t^3\tilde{K}^{\mp 2} +
t^{2\tilde{C} + 2} + t^{-2\tilde{C}-2}).$$

We have 
$$-\mathcal{X}^\pm\mathcal{X}^\mp(t^3 - t^{-3})^2 + t^{\mp 3}
\tilde{K}^2 + t^{\pm 3} \tilde{K}^{-2} = t^{2\tilde{C} + 2} +
t^{-2\tilde{C} - 2}.$$ \qed

Now we are in position to define the interpolating quantum group
$\U_{q,t}(G_2)$.

\begin{defi} $\U_{q,t}(G_2)$ is the algebra with generators $X_i^\pm$,
  $K_i^{\pm 1}$, $\tilde{K}_i^{\pm 1}$, $\eta$, $C$, $\tilde{C}$
  ($i = 1,2$) and relations ($i\neq j$)
$$\CC[K_i^{\pm},\tilde{K}_i^{\pm 1},\eta, C, \tilde{C}]_{i=1,2}
\quad \text{is commutative},$$
$$\mathcal{U}_1 = \langle X_1^\pm,K_1^{\pm 1},\tilde{K}_1^{\pm 1}
\rangle \simeq \U_{q,t}({}^LG_1),$$
$$\mathcal{U}_2 = \langle X_2^\pm, K_2^{\pm 1}, \tilde{K}_2^\pm 1,
\eta, C, \tilde{C}, \tilde{C}' \rangle \simeq \U_{q,t}(G_1),$$
$$K_i X_j^\pm = q^{\pm r_i C_{i,j}} X_j^\pm K_i\text{ , }\tilde{K}_i
X_j^\pm = t^{\pm r_iC_{i,j}/3}\tilde{K}_j,$$
$$[X_1^\pm , X_2^\mp] = 0.$$
\end{defi}

We define $\mathcal{X}_2^\pm$, $\mathcal{K}_2$ as for
$\U_{q,t}(G_1)$. Let $m\in\ZZ/2\ZZ$. From the above results we have
the following:

\begin{prop} The subalgebra of $\U_{q,1}(G_2)/(\tilde{K}_i = 1)$
  generated by $X_1^\pm$, $X_2^\pm/(q-q^{-1})$, $K_i^{\pm 1}$ is
  isomorphic to $\U_q(G_2)$.

The quotient by $\epsilon^{2\eta} = 1$, $K_1 = -1$, $K_2 = (-1)^m$,
$q^C = (-1)^m\epsilon$ of the subalgebra of $\U_{\epsilon,t}(G_2)$
generated by $X_1^\pm$, $\mathcal{X}_2^\pm$, $\tilde{K}_1^{\pm 1}$,
$\mathcal{K}_2^{\pm 1}$ is isomorphic to $\U_{-t}({}^LG_2)$.
\end{prop}

We have thus defined an interpolating quantum group $\U_{q,t}(\g)$ for
any simple Lie algebra $\g$. The same definition gives us such an
algebra for any symmetrizable Kac--Moody algebra $\g$ such that $r
\leq 3$. We just use the relations of the elementary interpolating
quantum groups in the same way as above. We conjecture that this
definition may also be generalized to all symmetrizable Kac--Moody
algebra $\g$.

\section{Representations of elementary interpolating quantum
  groups}\label{bun}

Representation theory of the elementary interpolating quantum
groups $\U_{q,t}(A_1)$, $\U_{q,t}(C_1)$, $\U_{q,t}({}^LG_1)$
is easily derived from the representation theory of the corresponding
standard quantum groups. So we need to consider only $\U_{q,t}(B_1)$
and $\U_{q,t}(G_1)$. For these algebras we will observe the simplest
examples of representations interpolating between finite-dimensional
representations of Langlands dual quantum groups of rank one. Even
though naively we have $^L B_1 = B_1$ and $^L G_1 = G_1$, the
resulting duality of representations is non-trivial.

\subsection{Simple finite-dimensional representations of
  $\U_{q,t}(B_1)$}\label{sfdr}

As in the representation theory of quantum groups, let us start with
Verma modules. We want to construct a Verma module $M(\la)$ with
respective highest eigenvalues of $(K, \tilde{K}, \eta, C,
\tilde{C})$ given by $\la =
(\la,\tilde{\la},E,A,\tilde{A})\in (\CC^*)^2\times \{\pm
1\}\times \CC^2$. We set:
$$M(\la) = \bigoplus_{n\geq 0}\CC (X^-)^n v_\la.$$
We have the obvious action 
$$X^-.((X^-)^n.v_{\lambda}) = (X^-)^{n+1}.v_{\lambda}\text{ , }C =
A\text{Id}\text{ , }\tilde{C} = \tilde{A}\text{Id}\text{ , }$$
$$K ((X^-)^n v_\lambda) = \la q^{-2n} ((X^-)^nv_{\lambda}) \text{
, } \tilde{K} ((X^-)^n v_\lambda) = \tilde{\la} t^{-n} ((X^-)^n
v_\lambda),$$ $$\eta ((X^-)^n v_{\lambda}) = (-1)^n X^+
(X^-)^n v_{\lambda}.$$ The point is to have a well-defined action of
$X^+$ such that $X^+v_{\lambda} = 0$. From the relation involving
$X^+X^-$, the action is uniquely defined. The relation involving
$X^-X^+$ evaluated on $v_\lambda$ imposes the following condition on
$\la$:

\begin{lem} The Verma module $M(\la)$ is non-trivial if and only if
$$q^{EA} t^{E\tilde{A}-\frac{1}{2}} \tilde{\la}^{-1} + q^{-EA}
t^{-E\tilde{A} + \frac{1}{2}} \tilde{\la} - q t^{-E\tilde{A} +
\frac{1}{2}} \tilde{\la} \la - q^{-1} t^{E\tilde{A}
-\frac{1}{2}} (\tilde{\la} \la)^{-1} = 0.$$
\end{lem}

Now we want to have a Verma module with a finite-dimensional
quotient. For the specializations defined above, we consider $p\geq
0$, $n = 2p$, $\la = q^n$, $\tilde{\la} = t^p$. So we obtain
the relation
$$q^{EA} t^{E\tilde{A}-\frac{1}{2}-p} + q^{-EA} t^{-E\tilde{A} +
  \frac{1}{2}+p} - q^{1+n} t^{-E\tilde{A} + \frac{1}{2} + p} -
q^{-1-n} t^{E\tilde{A} -\frac{1}{2} - p} = 0.$$
Thus, we have $EA = -(n + 1)$ or ($EA = (n+1)$ and $E\tilde{A} = p+1/2$). 
But to have the second specialization of Proposition \ref{inter}, 
we must have $Kq^{PC} = \epsilon$ at $q = \epsilon$, so 
$\epsilon^{n + EA} = \epsilon$. So we are in the second case $EA = n + 1$. 
Consider an $(n+1)$-dimensional
vector space
$$V_n = \CC v_0\oplus \CC v_1\oplus \cdots \oplus \CC v_n.$$ We denote
$n = 2p\in 2\ZZ$, $v_{-1} = v_{n+1} = 0$ and use the usual quantum number
notation $[m]_x = (x^m - x^{-m})/(x - x^{-1})$ for $m\in\ZZ$.

Let us consider operators on $V_n$ defined by $C = n+1$, $\tilde{C} =
p + 1/2$, and the following formulas:
$$X^+ v_{2j} = [n-2j + 1]_q v_{2j-1}\text{ , }X^+ v_{2j + 1} = [n-
  2j]_{qt}v_{2j},$$
$$X^- v_{2j} = [2j+1]_q v_{2j + 1} \text{ , }X^- v_{2j + 1} = [2j +
  2]_{qt} v_{2j + 2},$$
$$K.v_j = q^{n-2j} v_j\text{ , }\tilde{K} v_j = t^{p-j} v_j\text{ ,
}\eta.v_j = -j v_j.$$ The idea of this deformation is just to
replace the quantum box $[m]_q$ where $m$ is even by $[m]_{qt}$.

\begin{lem} 
The above formulas define an action of $\U_{q,t}(B_1)$ on $V_n$.
\end{lem}

\demo All relations are clear without computation, except for
relations (\ref{hun}). Let us check these relations.

$$(q-q^{-1})(qt - (qt)^{-1}) X^+X^-.v_{2j} = (q^{2j+1} -
q^{-2j-1})((qt)^{n-2j} - (qt)^{2j-n})v_{2j} $$ $$=(q^{n+1}t^{n-2j} +
q^{-n-1}t^{2j-n} - q^{n-4j-1}t^{n-2j} - q^{4j+1-n}t^{2j-n})v_{2j},$$

$$(q-q^{-1})(qt - (qt)^{-1}) X^-X^+ .v_{2j} = (q^{n-2j+1} -
q^{-n+2j-1})((qt)^{2j} - (qt)^{-2j}) v_{2j}
$$ $$= (q^{n+1}t^{2j} + q^{-n-1}t^{-2j} - q^{n+1-4j}t^{-2j}-
q^{4j-n-1}t^{2j}) v_{2j},$$

$$(q-q^{-1})(qt - (qt)^{-1}) X^+X^-.v_{2j+1} = ((qt)^{2j+2} -
(qt)^{-2j-2})(q^{n-2j-1} - q^{2j-n+1})v_{2j+1}
$$ $$= (q^{-n-1}t^{-2j-2} + q^{n+1}t^{2j+2} - q^{n-4j-3}t^{-2j-2} -
q^{4j+3-n}t^{2+2j})v_{2j+1},$$

$$(q-q^{-1})(qt - (qt)^{-1}) X^-X^+ .v_{2j+1} = ((qt)^{n-2j} -
(qt)^{-n+2j})(q^{2j+1} - q^{-2j-1}) v_{2j+1}
$$ $$= (q^{-n-1}t^{2j-n} + q^{n+1}t^{n-2j} - q^{n-1-4j}t^{n-2j}-
q^{4j-n+1}t^{2j-n}) v_{2j+1}.$$

The formulas are also satisfied at the limits as for $X^-X^+v_0$
and for $X^+X^-v_{2p}$ we get $q^{n+1} + q^{-n - 1} - q^{n+1} - q^{-n
  - 1} = 0$.
\qed

Note that $V_n' = \CC v_0\oplus \CC v_2\oplus\cdots \oplus \CC v_n$ is
stable for the action of $E^2$ and $F^2$.  We interpret this as a
Langlands duality of representations corresponding to $B_1$ and $C_1$
as follows: the first specialization of $V_n$ is the
$(n+1)$-dimensional simple representation of $\U_q(sl_2) = \U_q(B_1)$,
and by using the second specialization we obtain $V_n'$, the
$(p+1)$-dimensional representation of $\U_{t^2}(sl_2) = \U_{-t}(C_1) =
\U_{-t}(^L B_1)$ (at $q = \epsilon$, we have $K^2 = 1$ and 
$Kq^{PC} = \ep$ on $V_n$ as in Proposition \ref{inter}).

\begin{rem} If $n\notin 2\ZZ$, we can also define a representation 
of $\U_{q,t}(B_1)$ with the same formulas. 
Indeed the formulas are also satisfied at the limit: it is the same
for $X^-X^+v_0$ and for $X^+X^-v_{2j+1}$ where $n = 2j + 1$ we get 
$q^{-n-1}t^{-n-1} + q^{n+1}t^{n+1} - q^{-n-1}t^{-n-1} -
q^{n+1}t^{n+1} = 0$. But then we cannot use the second specialization
on $V_n'$ at $q = \epsilon$ as we have 
$Kq^{PC}v_0 = - v_0$ different than in Proposition \ref{inter}.
\end{rem}

\subsection{Representations of $\U_{q,t}(G_1)$}

Let $V_n$ be as in Section \ref{sfdr} where $n\in 3\ZZ$. Let us consider
operators defined by $C = n + 1$, $\tilde{C} = (n+1)/2$ and the
following formulas:
$$X^+ v_{3j} = (q-q^{-1})[n-3j + 1]_q v_{3j-1}\text{ , } X^- v_{3j} =
(q - q^{-1})[3j+1]_q v_{3j + 1},$$
$$X^+ v_{3j + 1} = (qt - (qt)^{-1})[n- 3j]_{qt}v_{3j}\text{ , }X^-
v_{3j + 1} = (q - q^{-1})[3j + 2]_q v_{3j + 2},$$
$$X^+ v_{3j+2} = (q - q^{-1})[n-3j-1]_qv_{3j+1}\text{ , }X^- v_{3j +
2} = (qt - (qt)^{-1})[3j + 3]_{qt} v_{3j + 2},$$
$$K.v_j = q^{n-2j} v_j\text{ , }\tilde{K} v_j = t^{n/2 - j} v_j\text{ ,
} \eta.v_{j} = -j v_j.$$ The idea of this deformation is just to
replace the quantum box $[m]_q$ where $m\equiv 0 [3]$ by
$[m]_{qt}$. This is analog to the deformation considered for $B_1$.

Note that we have in particular
$$P_+.v_{3j} = v_{3j}\text{ , }P_+.v_{3j+1} = 0\text{ , }P_+.v_{3j+2}
= -v_{3j+2},$$
$$P_-.v_{3j} = -v_{3j}\text{ , }P_-.v_{3j+1} = v_{3j+1}\text{ ,
}P_-.v_{3j+2} = 0.$$

\begin{lem} 
The above formulas define an action of $\U_{q,t}(G_1)$ on $V_n$.
\end{lem}

\demo All relations are clear without computation, except for
relations (\ref{hung}). Let us check these relations.

$$X^+X^-.v_{3j} = (q^{3j+1} - q^{-3j-1})((qt)^{n-3j} -
(qt)^{3j-n})v_{3j}
$$ $$=(q^{n+1}t^{n-3j} + q^{-n-1}t^{3j-n} - q^{n-6j-1}t^{n-3j} -
q^{6j+1-n}t^{3j-n})v_{3j},$$

$$X^+X^-.v_{3j+1} 
= (q^{3j+2} - q^{-3j-2})(q^{n-3j-1} - q^{3j-n+1})v_{3j+1}
$$ $$=(q^{n+1} + q^{-n-1} - q^{n-6j-3} - q^{6j+3-n})v_{3j+1},$$

$$X^+X^-.v_{3j+2} 
= ((qt)^{3j+3} - (qt)^{-3j-3})(q^{n-3j-2} - q^{3j-n+2})v_{3j+2}
$$ $$=(q^{n+1}t^{3j+3} + q^{-n-1}t^{-3j-3} - q^{n-6j-5}t^{-3j-3} -
q^{6j+5-n}t^{3j+3})v_{3j+2},$$

$$X^-X^+.v_{3j} 
= (q^{n-3j+1} - q^{3j-1-n})((qt)^{3j} - (qt)^{-3j})v_{3j}
$$ $$=(q^{n+1}t^{3j} + q^{-n-1}t^{-3j} - q^{n-6j+1}t^{-3j} -
q^{6j-1-n}t^{3j})v_{3j},$$

$$X^-X^+.v_{3j+1} 
= ((qt)^{n-3j} - (qt)^{3j-n})(q^{3j+1} - q^{-3j-1})v_{3j+1}
$$ $$=(q^{n+1}t^{n-3j} + q^{-n-1}t^{3j-n} - q^{n-6j-1}t^{n-3j} -
q^{6j+1-n}t^{3j-n})v_{3j+1},$$

$$X^-X^+.v_{3j+2} 
= (q^{n-3j-1} - q^{3j+1-n})(q^{3j+2} - q^{-3j-2})v_{3j+2}
$$ $$=(q^{n+1} + q^{-n-1} - q^{n-6j-3} - q^{6j+3-n})v_{3j+2}.$$
The formulas are also satisfied at the limits as for $X^-X^+v_0$
and also for $X^+X^-v_{3(n/3)}$ we get $q^{n+1} + q^{-n - 1} - q^{n+1}
- q^{-n - 1} = 0$.
\qed

Note that $V_n' = \CC v_0\oplus \CC v_3\oplus\cdots \oplus \CC v_n$ is
stable for the action of $E^3$ and $F^3$.  We interpret this as a
Langlands duality of representations of $G_1$: by using the first
specialization, $V_n$ becomes the $(n+1)$-dimensional simple
representation of $\U_q(sl_2) = \U_q(G_1)$, and by using the second
specialization we obtain $V_n'$, the $n/3 + 1$-dimensional representation
of $\U_{t^3}(sl_2) = \U_t({}^LG_1)$ (at $q = \epsilon$, we have $\epsilon^\eta = 1$,
$K = (-1)^{n/3}$ and $q^C = (-1)^{n/3} \ep$ on $V_n'$ as in Proposition \ref{speg}).

\begin{rem}
If $n\notin 3\ZZ$, we can also define a representation of
$\U_{q,t}(B_1)$ by the same formulas.  Indeed, the formulas are also
satisfied at the limit: it is the same for $X^-X^+v_0$, for
$X^+X^-v_{3j+1}$ where $n = 3j + 1$ we get $q^{n+1}+ q^{-n-1} -
q^{-n-1} - q^{n+1} = 0$, and for $X^+X^-v_{3j+2}$ where $n = 3j + 2$
we get $q^{n+1}t^{n+1}+ q^{-n-1}t^{-n-1} - q^{-n-1}t^{-n-1} -
q^{n+1}t^{n+1} = 0$.  But then we cannot use the second
specialization on $V_n'$ as at $q = \epsilon$ we have $Kq^C v_0 =
\epsilon^{2n+1} v_0$ different than in Proposition \ref{speg}.
\end{rem}

\subsection{Another interpretation of the duality}\label{disc}

In this section we discuss an interpretation of the Langlands dual of
Section \ref{sfdr} in terms of the structure of the algebra
$\U_{q,t}(B_1)$.

The duality of the simple finite dimensional representations of
$\U_{q,t}(B_1)$ in Section \ref{sfdr} in terms of characters is just
the elementary duality between the following polynomials:
$$y^{2n} + y^{2n-2} + \cdots + y^{2 - 2n} + y^{-2n}\leftrightarrow
(y^2)^n + (y^2)^{n-2} + \cdots + (y^2)^{2-n} + (y^2)^{-n}.$$ We have
seen that it corresponds to an interpolating representation. At the
level of characters, we can define a similar interpolation. Indeed let
us consider a map $\alpha(q,t)$ such that $\alpha(q,1) = 1$ and
$\alpha(\epsilon, t) = 0$. Such a map is naturally introduced in
\cite{fr, fr2} (we will also see an elementary way to introduce it
bellow):
$$\alpha(q,t) = (q + q^{-1})(qt - q^{-1}t^{-1})(q^2t -
q^{-2}t^{-1})^{-1}.$$ In the following it will just be denoted by
$\alpha$. Then the character
$$y^{2n} + \alpha y^{2n-2} + y^{2n-4} + \alpha y^{2n-6} + \cdots +
y^{4 - 2n}+ \alpha y^{2 - 2n}+ y^{-2n}$$ interpolates between the two
formulas.

The character of a fundamental representation is $y + y^{-1}$. This
corresponds to the decomposition of the Casimir central elements
$$Cas(q) = q^C+q^{-C}, \qquad Cas(t) = t^{\tilde{c}+
1/2}+t^{-\tilde{c} - 1/2}.$$ The Casimir element of the second
specialization is $(t^{2(\tilde{c}+1)}+t^{-2(\tilde{c}+1)})$, so we have the
following picture:

$$\begin{CD}
\U_{q,t}(B_1)  @>{ t \rightarrow 1}>>  \U_q(B_1)\otimes_{Z_q(B_1)}
\CC[q^C,q^{-C}]\\
@VV{q\rightarrow \epsilon}V     @VV{\text{Duality}}V\\
\U_{\epsilon ,t}(B_1)  @<<{\text{Injection}}<
\U_{-t}(C_1)\otimes_{Z_{-t}(C_1)} \CC[t^{\tilde{c}},t^{-\tilde{c}}]
\end{CD}.$$
Note that the tensor product $\U_q(B_1)\otimes_{Z_q(B_1)}
\CC[q^{C},q^{-C}]$ is a quantum analogue (for $\g=B_1$) of the algebras
$$
U(\g) \otimes_{Z(U(\g))} U({\mathfrak h}),
$$
where ${\mathfrak h} \subset \g$ is the Cartan subalgebra,
considered by Gelfand and Kirillov \cite{GK}.

The decomposition of the Casimir element and the character formulas
are closely related. This can be put a step forward by having a
similar interpretation of the interpolating character in the spirit of
the constructions of \cite{fr, fr2} in the affine case. Indeed, we
have the central element $Y + Y^{-1}$ in $\U_{q,t}(B_1)$ where $Y =
q^{PC}t^{\tilde{c}+ 1/2}$, which interpolates between $Cas(q)$ and
$Cas(t)$. (Note that we have $Y^{\pm 1} X^+ = X^+Y^{\mp 1}$.)
 
We define the completed algebra $\tilde{\U}_{q,t}(B_1)$ as the algebra
containing $\U_{q,t}(B_1)$ with additional elements $W^{\pm 1}$ such
that $W^{\pm 1}X^+ = X^+W^{\mp 1}$ and $YW =\alpha WY$. Note that
because of the relation of the algebra, we cannot require it to
commute with $X^-$. Note also that we have $Y^{-1}W^{-1} = \alpha
W^{-1}Y^{-1}$ which is compatible with the commutation relations with
$X^+$.

Let us explain how such a commutation relation $YW =\alpha WY$ can be
obtained naturally in the spirit of \cite{fr, fr2}. We recall that the
variables for the $q$-characters (affine version of characters) are
materialized as formal power series in generators of the Heisenberg
subalgebra of the level $0$ quantum affine algebra. The
$q,t$-analogues of these variables (which are the building blocks for
the generating series of the deformed ${\mc W}$-algebra ${\mc
W}_{q,t}(\g)$) are, in turn, related to a non-commutative Heisenberg
algebra. In the $B_1$-case this Heisenberg algebra has generators
$h[n]$ ($n\in\ZZ$) such that $h[0]$ is central and for $n,m\neq 0$ we
have
$$[h[n],h[m]] = \delta_{n,-m} \frac{(q^{n} - q^{-n}) (t^n -
  t^{-n})}{n}.$$

As the $0$th mode $h[0]$ is central, it is not clear how to obtain the
commutation relations as considered above. But the finite type can
also be seen as a limit of the affine type case in the following
sense. Let $$h_\pm(z) = \text{exp}\left(\sum_{\pm m\geq 0}h[m] z^{-m}\right).$$
We have
$$h_+(zq^3t)h_-(w) = h_-(w)h_+(zq^3t)\text{exp}\left(\sum_{m > 0}
\frac{(q^{-2m} - q^{-4m}) (1 - t^{-2m})}{m} (wz^{-1})^m\right)$$
$$= h_-(w)h_+(zq^3t) \frac{(1 - q^{-4}wz^{-1})(1 -
t^{-2}q^{-2}wz^{-1})}{(1 - q^{-4}t^{-2}wz^{-1})(1 - q^{-2}wz^{-1})}.$$
A priori, we cannot directly specialize at $z = w = 1$. But if we
forget the intermediate formulas, everything makes sense for this
specialization, and for $Y = h_+(q^3t)$, $W = h_-(1)$ we obtain
$$YW =  \frac{(1 - q^{-4})(1 - t^{-2}q^{-2})}{(1 - q^{-4}t^{-2})(1 -
  q^{-2})}WY = \alpha WY.$$ 

To give a precise meaning to this specialization, we consider an
additional formal parameter $u$ and replace $w$, $z$, respectively, by
$wu$, $zu^{-1}$. We get formal power series in $u^{-1}$. So we can set
$z = w = 1$ and for $Y(u) = h_+(u^{-1}q^3t)$ and $W(u) = h_-(u)$ we
obtain
$$Y(u)W(u) = \frac{(1 - q^{-4}u^2) (1 - t^{-2}u^2)}{(1 - q^{-2}u^2)(1
  - q^{-4}t^{-2}u^2)}W(u)Y(u).$$
Now we can specialize from the affine type to the finite type by
considering $Y = Y(1)$, $W = W(1)$, and we get $YW = \alpha WY$ as
explained above.

We have a notion of normal ordering $:M:$ for monomials $M$ in
$Y^{\pm 1}, W^{\pm 1}$, where we put the $Y^{\pm 1}$ on the left and
the $W^{\pm 1}$ on the right. Then we have
$$(:YW + Y^{-1}W^{-1}:)^2
= \alpha^{-1}:(YW)^2: + 2 \alpha + \alpha^{-1} :(Y^{-1}W^{-1})^2:$$
$$ = \alpha^{-1} (:(YW)^2: + \alpha^2 + :(Y^{-1}W^{-1})^2:) +
\alpha.$$
In particular, the formula corresponding to the $3$-dimensional simple
representation appears naturally as
$$:(YW)^2: + \alpha^2 + :(Y^{-1}W^{-1})^2:$$ Note that this formula
commutes with $X^+$, as does $:YW + Y^{-1}W^{-1}:$, which has the same
property. This can be interpreted as an analog of the invariance of
the usual characters for the Weyl group action or of the symmetry
property of $q$-characters for the screening operators (see
\cite{fr2,fr3}).

It would be desirable to have a similar interpretation of the duality
of characters for general quantum groups.

\section{More general interpolating representations}    \label{int reps}

By an {\em interpolating representation} we understand a
representation of the interpolating quantum group $\U_{q,t}(\g)$ which
gives by specialization representations of the Langlands dual quantum
groups. We have seen in the Section \ref{bun} that interpolating
representations exist for elementary interpolating quantum groups. In
this section we give additional examples for non-elementary
interpolating quantum groups. We believe that any irreducible
representation $L(\la)$ of $U_q(\g)$ (equivalently, of $\U_q(\g)$)
with $\la \in P'$ may be $t$-deformed, in an essentially unique way,
to a representation of $\U_{q,t}(\g)$ in such a way that its
specialization at $q=\ep$ gives a representation of
$U_{-t}({}^L\g)$ whose character is $\Pi(\chi(\lambda))$.

We start with a simple finite-dimensional representation $V$ of
$\U_q(B_2)$ with highest weight which has an even multiplicity for the
node $2$. We want to "deform" the $\U_q(B_2)$-module structure on
$V$. All weights of $V$ have even multiplicities for the node $2$. For
$v\in V$ of weight $m_1 \om_1 + 2 m_2\om_2$, we set
$$K_2 v = q^{2 m_2}v\text{ , }\tilde{K}_2 v = t^{m_2}v\text{ , }K_1 v
= q^{2m_1}v\text{ , }\tilde{K}_1 v = t^{m_1}v.$$ The deformation will
be necessarily semi-simple for $\U_1\simeq \U_{q^2t}(A_1)$, but
moreover we require that it is semi-simple for the action of $\U_2 =
\U_{q,t}(B_1)$ with simple submodules isomorphic to the
representations constructed in Section \ref{bun}. The actions of $C$
and $\tilde{C}$ are uniquely determined from the action of $X^+_2$ and
$X^-_2$ as it suffices to know the decomposition in simple modules for
$\U_2$. So the non-trivial point is to deform the action of the
$X^+_i$, $X^-_i$.

We will consider $3$ examples of interpolating representations of
$\U_{q,t}(B_2)$. At the level of crystals, they correspond to the
examples studied in Section \ref{charstate}. The first one is the most
simple example where the duality occurs. In the second example we have
a multiplicity in the character and we can see that the relations
between $C$ and the $X^\pm_1$ cannot be written a simple way. In the
third example we observe that different $t$-deformations of the Serre
relations arise in the interpolating representations.

\mk

{\bf Example 1.} Let $V = L(\om_1)$ be the fundamental representation
of $\U_q(B_2)$ of dimension $5$ which corresponds by duality to the
representation of $\U_{-t}(C_2)$ whose highest component is the
fundamental representation of dimension $4$. Its character is $y_1 +
y_2^2y_1^{-1} + 1 + y_1y_2^{-2} + y_1^{-1}$, and all weight spaces are
of dimension $1$.

We consider a basis $(v_l)_{1\leq l\leq 5}$ of $V$ such that $v_1$ is
a highest weight vector,
$$v_2 = X^-_1 v_1\text{ , }v_3 = X^-_2v_2\text{ , }v_4 = X^-_2
v_3/[2]_q\text{ , }v_5 = X^-_1 v_4.$$ In this basis the action of the
$X^\pm_i$ has matrix coefficients $0$, $1$ or $[2]_q$. We
deform the action by replacing the $[2]_q$ by $[2]_{qt}$, that is to
say we only deform $X^{\pm}_2v_3 = [2]_{qt} v_{3\mp 1}$. The
decomposition in simple modules for $\U_1$ and $\U_2$ is clear and
coincides with the case $t = 1$.

\mk

{\bf Example 2.} Let $V = L(2\om_2)$ representation of $\U_q(B_2)$ of
dimension $10$ which corresponds by duality to the representation of
$\U_{-t}(C_2)$ whose highest component is the fundamental
representation of dimension $5$.  Its character is
$y_2^2 + y_1 + y_1^2y_2^{-1} + y_2^2y_1^{-1} + 2 . 1 + y_2^2y_1^{-2} +
y_1y_2^{-1} + y_1^{-1} + y_2^{-2}$. There is a multiplicity $2$ for
the weight $1$.

We consider a basis $(v_l)_{1\leq l\leq 10}$ of $V$ such that $v_1$ is
a highest weight vector,
$$v_2 = X^-_2 v_1\text{ , }v_3 = X^-_2v_2/[2]_q\text{ , }v_4 = X^-_1
v_3\text{ , }v_5 = X^-_1 v_4/[2]_{q^2}\text{ , }v_6 = X^-_2 v_5,$$
$$v_7 = X^-_2 v_6/[2]_q\text{ , }v_8 = X^-_1 v_2\text{ , }v_9 = X^-_2
v_8\text{ , }v_{10} = X^-_2 v_9/[2]_q.$$ In this basis the action of
the $X^\pm_i$ have matrix coefficients $0$, $1$, $[2]_q$ or
$[2]_{q^2}$. We deform the action by replacing these coefficients,
respectively, by $0$, $1$, $[2]_{qt}$, $[2]_{q^2t}$. That is to say we
only deform
$$X^\pm_1 v_4 = [2]_{q^2t} v_{4\mp 1}\text{ , }X^\pm_2v_6 =
[2]_{qt}v_{6\mp 1}\text{ , }X^\pm_1v_9 = [2]_{qt} v_{4\mp 1}.$$ The
decomposition in simple modules for $\U_1$ and $\U_2$ is clear and
coincides with the case $t = 1$ except for the trivial submodules of
$\U_2$ and $\U_1$ which are, respectively,
$$\CC(v_9 - [2]_{qt}v_4)\text{ and }\CC([2]_{qt}v_4 -
[2]_{q^2t}v_9).$$ Note that a priori we cannot expect to have simple
relations between the $C$ and the $X^\pm_1$ as $v_4$ is not an
eigenvector of $C$.  \mk

{\bf Example 3.} Let $V = L(2\om_1)$, an irreducible representation of
$\U_q(B_2)$ of dimension $14$, which corresponds by duality to a
representation of $\U_{-t}(C_2)$ whose highest component is of
dimension $10$. Its character is $y_1^2 + y_2^2 + y_2^4y_1^{-2} + y_1
+ y_2^2y_1^{-1} + y_1^2y_2^{-2} + 2 . 1 + y_1y_2^{-2} + y_2^2 y_1^{-2}
+ y_1^2 y_2^{-4} + y_1^{-1} + y_2^{-2} + y_1^{-2} $, and there is a
multiplicity $2$ for the weight $1$.

We consider a basis $(v_l)_{1\leq l\leq 14}$ of $V$ such that $v_1$ is
a highest weight vector,
$$v_2 = X^-_1 v_1\text{ , }v_3 = X^-_1v_2/[2]_{q^2}\text{ , }v_4 =
X^-_2 v_3\text{ , }v_5 = X^-_2 v_4/[2]_q,$$
$$v_6 = X^-_2 v_5/[3]_q\text{ , }v_7 = X^-_2 v_6/[4]_q\text{ , }v_8 =
X^-_1 v_7\text{ , }v_9 = X^-_1 v_8/[2]_{q^2}\text{ , }v_{10} = X^-_2
v_2,$$
$$v_{11} = X^-_2 v_{10}/[2]_q\text{ , }v_{12} = X^-_1 v_{11}\text{ ,
}v_{13} = X^-_1 v_{12}/[2]_{q^2}\text{ , }v_{14} = X^-_1 v_6.$$ In
this basis the action of the $X^\pm_i$ have matrix coefficients
$0$, $1$, $[2]_q$, $[3]_q$, $[4]_q$, $[2]_{q^2}$,
$([2]_q[2]_{q^2}/[4]_q)^{\pm 1}$, $[4]_q/[2]_{q^2}$. We deform the
action by replacing these coefficients, respectively, by $0$, $1$,
$[2]_{qt}$, $[3]_q$, $[4]_{qt}$, $[2]_{q^2t}$,
$([2]_{qt}[2]_{q^2t}/[4]_{qt})^{\pm 1}$, $[4]_{qt}/[2]_{q^2t}$. The
decomposition into simple modules for $\U_1$ and $\U_2$ is clear and
coincides with the case $t = 1$ except for the trivial submodule of
$\U_2$ and $\U_1$ which are respectively
$$\CC([3]_q[4]_{qt} v_{12} - [2]_{qt}[2]_{q^2t}v_5)\text{ and
}\CC([2]_{qt}[2]_{q^2t}^2v_{12} - [4]_{qt}v_5).$$

In this example we can observe non-trivial $t$-deformations of the
Serre relations, but different relations are satisfied on different
vectors in the representation! Indeed, we have
$$(X^-_2{X^-_1}^2 - (q^2 t^2 + q^{-2}t^{-2})X^-_1X^-_2X^-_1 +
{X^-_1}^2 X^-_2)v_1 = 0,$$
$$(X^-_2{X^-_1}^2 - (q^2t + q^{-2}t^{-1})X^-_1X^-_2X^-_1 +
{X^-_1}^2X^-_2)v_{11} = 0.$$ This implies that if we impose any
$t$-deformation of the Serre relation
$$(X^-_2{X^-_1}^2 - (q^2 + q^{-2})X^-_1X^-_2X^-_1 +
{X^-_1}^2X^-_2) = 0$$
in the algebra $\U_{q,t}(B_2)$, then the resulting algebra will not
act on the $\U_{q,t}(B_2)$-module that we have just constructed.
Indeed, this relation will be different from the relation satisfied on
at least one of the vectors, $v_1$ and $v_{11}$ (as written
above). The difference of the two relations would give us a
multiple of the monomial $X^-_1X^-_2X^-_1$, which should then
have to annihilate this vector. But neither vector is annihilated by
this monomial: we have
$$
(X^-_1X^-_2X^-_1) v_1    = v_4 [2]_{qt}[2]_{q^2t}/[4]_{qt},
$$
$$
(X^-_1X^-_2X^-_1) v_{11} = v_{14} [2]_{qt}[2]_{q^2t}/[4]_{qt}.
$$
Moreover, one can show that the structure of $\U_{q,t}(\g)$-module on
$V$ described above is unique (the same is also true for the modules
in Examples 1 and 2).

\bigskip

Now let us explain how we constructed the above interpolating
representations.  Let $V$ be a simple finite-dimensional
representation of $\U_q(B_2)$ as above of highest weight $\lambda =
m_1\om_1 + 2m_2\om_2$. We have a decomposition in weight spaces $V =
\bigoplus_{\mu\leq \lambda} V_\mu$. Let $V_r = \bigoplus_{\mu =
\lambda - \alpha_{i_1} - \cdots - \alpha_{i_r}} V_{\mu}$ and $V_{\leq
R} = \bigoplus_{r\leq R} V_r$. We have $X^-_1V_R + X^-_2V_R = V_{R+1}$
and $V_\lambda = X^-_1 V_{\lambda - \alpha_1} + X^-_2 V_{\lambda -
\alpha_2}$.

We define on $V$ the action of the $K_i(t)$, $\tilde{K}_i(t)$ as
explained above.

We define by induction on $r\geq 0$ the deformed actions
$$X^+_1(t),X^+_2(t) : V_{r+1} \rightarrow V_r\text{ and
}X^-_1(t),X^-_2(t) : V_r\rightarrow V_{r+1},$$ satisfying the
following properties:

(i) $X^+_1(1) = X^+_1$, $X^+_2(1) = X^+_2$, $X^-_1(1) = X^-_1$,
$X^-_2(1) = X^-_2$,

(ii) the action is compatible with the restrictions to $\U_1$ and
$\U_2$,

(iii) $[X^+_1(t),X^-_2(t)] = 0$,

(iv) $[X^+_2(t),X^-_1(t)] = 0$.

To start with we set $X^+_1(t) = 0$, $X^+_2(t) = 0$ on $V_0$. 

Suppose that the deformed action is defined for $r\geq 0$. Let
$V_{\lambda}\subset V_{r+1}$. We want to define the deformed actions
$$\xymatrix{V_{\lambda +\alpha_1}\ar[dr]^{X^-_1(t)}& & V_{\lambda +
    \alpha_2}\ar[dl]^{X^-_2(t)}& & V_{\lambda + \alpha_1} \\ &
    V_\lambda & & V_\lambda \ar[ul]^{X^+_2(t)} \ar[ur]^{X^+_1(t)}& }$$
    By using the condition (ii) for $\U_2$, we can first define the
    action of $X^+_2(t)$ and $X^-_2(t)$. This gives in particular a
    decomposition $V_\lambda = V^{(2)}_\lambda\oplus
    \tilde{V}^{(2)}_\lambda$ where $V^{(2)}_\lambda = {X^-_2}(t)
    (V_{\lambda + \alpha_2})$ and $\tilde{V}^{(2)}_\lambda =
    \text{Ker}(X^+_2(t))\cap V_\lambda$. The condition (iii) on
    $V_{\lambda + \alpha_2}$ gives $\phi(t) :
    V^{(2)}_\lambda\rightarrow V_{\lambda + \alpha_1}$. The condition
    (iv) on $V_{\lambda + \alpha_1}$ gives $\psi(t) : V_{\lambda +
    \alpha_1}\rightarrow V^{(2)}_\lambda$. So it suffices to define
    $X^\pm_1(t)$ such that $X^+_1(t) = \phi(t)$ on $V^{(2)}_\lambda$,
    $\Pi\circ X^-_1(t) = \Psi(t)$ where $\Pi$ is the projection on
    $V^{(2)}_\lambda$ along $\tilde{V}^{(2)}_\lambda$, and
    $X^+_1(t)X^-_1(t) = R(t)$ given by condition (ii) for $\U_1$. In a
    matrix form we have $X^+_1(t) = \begin{pmatrix}\phi(t) &
    A(t)\end{pmatrix}$, $X^-_1(t) =
    \begin{pmatrix}\Psi(t)\\B(t)\end{pmatrix}$, $X^-_1(t)X^+_1(t) =
    \phi(t)\psi(t) + A(t)B(t)$.

So it suffices to prove that 
$$\text{rk}(R(t) - \phi(t)\psi(t)) \leq
\text{dim}(\tilde{V}^{(2)}_\lambda).$$ We call this the {\em
compatibility condition}. In the examples studied above this condition
is satisfied, and that is why the interpolating representations do
exist. We conjecture that it is satisfied in general and we have
the following

\begin{conj}\label{modconj} For any $\la \in P'$ there exists a unique irreducible
  representation $L_{q,t}(\la)$ of $\U_{q,t}(\g)$ whose specialization
  at $t = 1$, viewed as a $\U_q(g)$-module, is the irreducible module
  $L(\la)$ and specialization at $q = \ep$, viewed as a
  $\U_{-t}({}^L\g)$-module, contains a module of character
  $\Pi(\chi(\lambda))$.
\end{conj}

\section{Conjectures on the Langlands duality for quantum groups}\label{conj}

In this Section we conjecture stronger statements on the duality for
characters and crystals which we prove for simply-laced $\g$ with $r =
2$ and for $B_2$. The proof of these conjectures and the computation of
the corresponding Langlands duality branching rules is a program for
further study for this Langlands duality.

\subsection{A positivity conjecture} 

It is easy to compute the Langlands duality branching rules for the
examples of Section \ref{charstate}.

For $\Glie = B_2$:
$$\Pi(\chi(\om_1)) = \chi^L(\check\om_1) + \chi^L(0)\text{ , }
\Pi(\chi(2\om_2)) =
\chi^L(\check\om_2) + \chi^L(\check\om_1) + \chi^L(0),$$
$$\Pi(\chi(2\om_1)) = \chi^L(2\check\om_1) +
\chi^L(\check\om_1),$$
$$\Pi(\chi(\om_1 + 2\om_2)) 
= \chi^L(\check\om_1 + \check\om_2) + \chi^L(2 \check\om_1) +
\chi^L(2\check\om_2) + \chi^L(\check\om_1).
$$
and for $\Glie=G_2$:
$$\Pi(\chi(\om_1)) = \chi^L(\check\om_1) + \chi^L(0)\text{ , }
\Pi(\chi(\om_2)) = \chi^L(\check\om_2) + 2\chi^L(\check\om_1) +
\chi^L(0).$$
So it is natural to give a purely classical analog to Conjecture \ref{modconj}:
\begin{conj}\label{lchar} For any $\lambda\in P^+\cap P'$,
  $\Pi(\chi(\lambda))$ is the character of an $^L\g$-module.
\end{conj}
This Conjecture means that the virtual representation of Proposition
\ref{debut} is an actual representation, that is, the Langlands
duality branching rules are positive:
$$
\Pi(\chi(\lambda)) = \sum_{\check\mu \in P^{L,+}} m_{\check\mu}
\chi^L(\check\mu), \qquad m_{\check\mu} \in \Z_+.
$$
We will prove the conjecture in several cases, but first we prove that
in general certain Langlands duality branching rules are positive. We
use the partial ordering on $P'$ viewed as the ${}^L\g$ weight
lattice.

\begin{prop} Let $\check\mu_0$ maximal in $\{\check\mu\in
  P^{L,+}|m_{\check\mu}\neq 0, \check\mu\neq\Pi(\lambda)\}$.  Then
$m_{\check\mu_0} > 0$.
\end{prop}

\demo By Theorem \ref{incchar} the coefficient of $\check\mu_0$ in
$\Pi(\chi(\lambda))$ is larger than in $\chi^L(\Pi(\lambda))$. But the
only terms which can contribute to this multiplicity are by hypothesis
$\chi^L(\Pi(\lambda))$ and $\chi^L(\check\mu_0)$. This implies the
result.  \qed

This includes all coefficients in the examples at the beginning of
this section. Now let us consider a statement analogous to Conjecture
\ref{lchar} in terms of crystals.

\subsection{Structure of the crystal $\tilde{\mathcal{B}}(\lambda)$} 

For $\lambda\in P^+\cap P'$ let
$$\tilde{\mathcal{B}}(\lambda) = \{v\in\mathcal{B}(\lambda)|
\mathrm{wt}(m')\in P'\}.$$ Note that for $\Glie$ of type $B_\ell$ we
have $\tilde{\mathcal{B}}(\lambda) = \mathcal{B}(\lambda)$.

$\tilde{\mathcal{B}}(\lambda)\sqcup \{0\}$ is stable under the action
of the operators $e_i^L,f_i^L$. We define new maps
$\epsilon_i^L,\phi_i^L$ on $\tilde{\mathcal{B}}(\lambda)$
corresponding to the operators $e_i^L,f_i^L$ (in general they do not
coincide with the original $\epsilon_i,\phi_i$, as we may have
$e_i^L(m') = 0$ but $e_i^L(m')\neq 0$, see the example in Section
\ref{excrys}). We get an abstract $^L\g$-crystal $(\tilde{\mathcal{B}}
(\lambda),e_i^L,f_i^L,\epsilon_i^L,\phi_i^L,
\mathrm{wt}^L)$. Moreover, this crystal is semi-normal, that is to say
that for each $i\in I$, it is as a $^L\g_{\{i\}}$-crystal the crystal
of a $^L\g_{\{i\}}$-module (this is a direct consequence of the
semi-normal property of $\mathcal{B}(\lambda)$).  Here for $J\subset
I$ we denote by $\Glie_J$ the sub Lie algebra of Cartan matrix
$(C_{i,j})_{i,j\in J}$.

Theorem \ref{dcrys} means that the connected component of
$\tilde{\mathcal{B}}(\lambda)$ containing the highest weight vector is
normal, that it to say is the crystal of a $^L\g$-module. In this
section we study the structure of the whole crystal
$\tilde{\mathcal{B}}(\lambda)$.

In all examples of Section \ref{excrys} the crystal is normal. In
particular, we have obtained the following:

\begin{prop}\label{pfund} For all pseudo
fundamental representations of a rank $2$ Lie algebra,
$\tilde{\mathcal{B}}(\lambda)$ is normal.
\end{prop}

So we could expect naively that $\tilde{\mathcal{B}}(\lambda)$
is normal. This statement is not true in
general. For example in type $B_2$ consider $\lambda = \om_1 +
2\om_2$.

We have seen that in terms of characters $\Pi(\chi(\om_1 + 2\om_2))$ has $4$ simple
constituents.
But $\mathcal{B}(\lambda) =
\mathcal{M}(Y_{1,0}Y_{2,1}^2) = \tilde{\mathcal{B}}(\lambda)$ has $3$
connected component as a $^L\g$-crystal.

The first connected component is isomorphic to
$\mathcal{B}^L(\check\om_1 + \check\om_2)$ ($16$ terms):
\\$\{1_02_1^2, 1_2^{-1}2_1^4, 1_01_2^22_3^{-2}, 2_1^22_3^{-2}1_2,
1_01_21_4^{-1}, 2_1^21_4^{-1}, 1_01_4^{-2}2_3^2, 2_3^{-4} 1_2^3, \\
1_2^{-1}1_4^{-2}2_1^22_3^2, 2_3^{-2}1_2^21_4^{-1}, 1_02_5^{-2},
1_21_4^{-2}, 1_2^{-1}2_1^22_5^{-2}, 1_4^{-3}2_3^2,
2_3^{-2}2_5^{-2}1_2, 1_4^{-1}2_5^{-2}\}$.

The second connected component is isomorphic to $\mathcal{B}^L(\om_1)$
($4$ terms): \\$\{1_02_12_5^{-1}, 1_2^{-1}2_1^32_5^{-1},
2_12_3^{-2}2_5^{-1}1_2, 1_4^{-1}2_12_5^{-1}\}$.

The third connected component is ($15$ terms):
\\$\{1_01_22_12_3^{-1}, 2_1^32_3^{-1}, 1_01_4^{-1}2_12_3,
1_2^{-1}1_4^{-1}2_1^32_3, 2_1 2_3^{-3}1_2^2, 1_01_22_3^{-1}2_5^{-1},
2_12_3^{-1}1_21_4^{-1}, 2_1^22_3^{-1}2_5^{-1},
\\1_01_4^{-1}2_32_5^{-1}, 1_2^{-1}1_4^{-1}2_1^22_32_5^{-1},
1_4^{-2}2_12_3, 1_2^22_3^{-3}2_5^{-1}, 2_12_3^{-1}2_5^{-2},
1_21_4^{-1}2_3^{-1}2_5^{-1}, 1_4^{-2}2_32_5^{-1}\}$.

Although the third component has two highest weight elements $u =
1_01_22_12_3^{-1}$ and $v = 2_1^32_3^{-1}$ it is not connected as
$e_2^Le_1^Le_1^Le_2^L v =
1_21_4^{-1}2_3^{-1}2_5^{-1} =
e_1^Le_2^Le_2^Le_1^Le_1^Lu$. But
its character is $\chi^L(\check\om_1) + \chi^L(\check\om_2)$.

In fact, it suffices to modify slightly the
crystal structure of the third component to get a normal crystal. 
Indeed, without changing the $\mathrm{wt}^L,
\epsilon_i^L,\phi_i^L$, we just replace
$$\begin{cases}e_2^L 2_1^32_3^{-1} &= 2_12_3^{-3}1_2^2,
\\e_2^L 1_01_4^{-1}2_12_3 &= 1_01_22_3^{-1}2_5^{-1},\end{cases}
\text{ by } 
\begin{cases}e_2^L 2_1^32_3^{-1} &= 1_01_22_3^{-1}2_5^{-1},
\\e_2^L 1_01_4^{-1}2_12_3 &= 2_12_3^{-3}1_2^2.\end{cases}$$ In
other words we have defined a bijection $\Psi :
\tilde{\mathcal{B}}(\lambda)\rightarrow \mathcal{B}^L$, where
$\mathcal{B}^L$ is normal, satisfying
$(\mathrm{wt}^L,\phi_i^L,\epsilon_i^L)\Psi = (\mathrm{wt}^L, \phi_i^L,
\epsilon_i^L)$ and $\Psi f_1^L = f_1^L\Psi$,
$\Psi e_1^L = e_1^L\Psi$. But $\Psi$ is not a morphism
of crystal as $\Psi f_2^L \neq f_2^L\Psi$ and
$\Psi e_2^L \neq e_2^L\Psi$ (see the picture bellow).

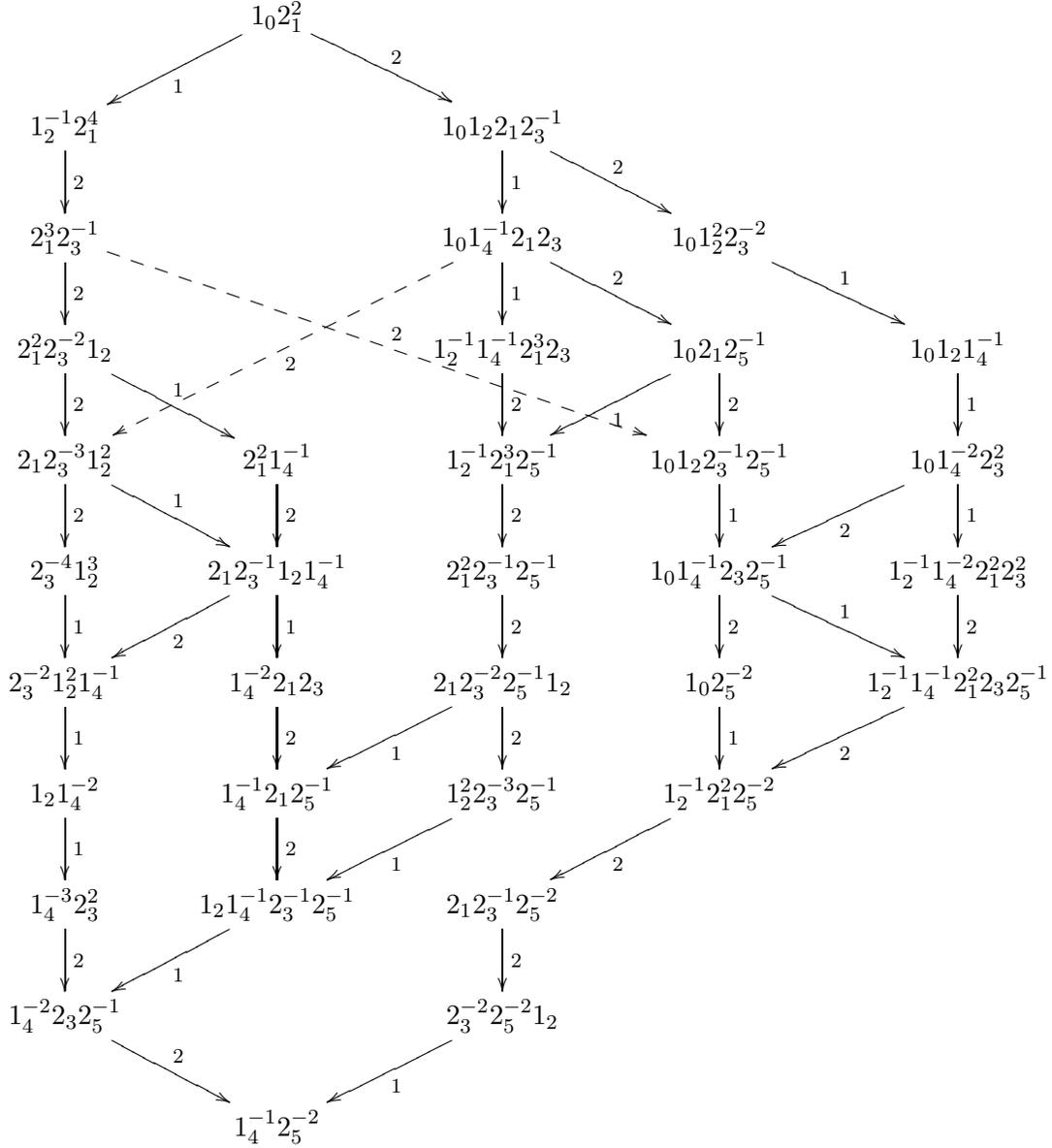
\xymatrix{ &1_02_1^2\ar^1[dl]\ar^2[dr]&&&
\\1_2^{-1}2_1^4\ar^2[d]&&1_01_22_12_3^{-1}\ar[d]^1\ar[dr]^2&&
\\2_1^32_3^{-1}\ar[d]^2\ar@{-->}[ddrrr]^2&&1_01_4^{-1}2_12_3\ar@{-->}[ddll]^2\ar[d]^1\ar[dr]^2&1_01_2^22_3^{-2}\ar[dr]^1&
\\2_1^22_3^{-2}1_2\ar[d]^2\ar[dr]^1&&1_2^{-1}1_4^{-1}2_1^32_3\ar[d]^2&1_02_12_5^{-1}\ar[dl]^1\ar[d]^2&1_01_21_4^{-1}\ar[d]^1
\\2_1
2_3^{-3}1_2^2\ar[dr]^1\ar[d]^2&2_1^21_4^{-1}\ar[d]^2&1_2^{-1}2_1^32_5^{-1}\ar[d]^2&1_01_22_3^{-1}2_5^{-1}\ar[d]^1&1_01_4^{-2}2_3^2\ar[d]^1\ar[dl]^2
\\2_3^{-4}1_2^3\ar[d]^1&2_12_3^{-1}1_21_4^{-1}\ar[dl]^2\ar[d]^1&2_1^22_3^{-1}2_5^{-1}\ar[d]^2&1_01_4^{-1}2_32_5^{-1}\ar[dr]^1\ar[d]^2&1_2^{-1}1_4^{-2}2_1^22_3^2\ar[d]^2
\\2_3^{-2}1_2^21_4^{-1}\ar[d]^1&1_4^{-2}2_12_3\ar[d]^2&2_12_3^{-2}2_5^{-1}1_2\ar[d]^2\ar[dl]^1&1_02_5^{-2}\ar[d]^1&1_2^{-1}1_4^{-1}2_1^22_32_5^{-1}\ar[dl]^2
\\1_21_4^{-2}\ar[d]^1& 1_4^{-1}2_12_5^{-1}\ar[d]^2&
1_2^22_3^{-3}2_5^{-1}\ar[dl]^1&1_2^{-1}2_1^22_5^{-2}\ar[dl]^2&
\\1_4^{-3}2_3^2\ar[d]^2&1_21_4^{-1}2_3^{-1}2_5^{-1}\ar[dl]^1&2_12_3^{-1}2_5^{-2}\ar[d]^2&&
\\1_4^{-2}2_32_5^{-1}\ar[dr]^2&&2_3^{-2}2_5^{-2}1_2\ar[dl]^1&&
\\&1_4^{-1}2_5^{-2}&&&}

\begin{conj}\label{twiscrys} For $\lambda\in P'$, there is a bijection
  $\Psi : \tilde{\mathcal{B}}(\lambda)\rightarrow \mathcal{B}^L$ to
  a normal ${}^L\g$-crystal $\mathcal{B}^L$ satisfying
  $\mathrm{wt}^L \Psi = \mathrm{wt}^L$ and
  $(\phi_i^L,\epsilon_i^L)\Psi = (\phi_i^L, \epsilon_i^L)$ for any
  $i\in I$.
\end{conj}

\noindent This means that, by changing the maps $e_i^L$,
$f_i^L$, respectively, by $\Psi^{-1}e_i^L\Psi$,
$\Psi^{-1}f_i^L\Psi$, we get a normal crystal.

\bigskip

Conjecture \ref{twiscrys} implies Conjecture \ref{lchar} as we have
$$\Pi(\chi(\lambda)) = \sum_{v\in\tilde{\mathcal{B}}(\lambda)}
\mathrm{wt}^L(v).$$

First, we look at the case of the Lie algebra of rank $1$.  For $r =
1$ the result is clear as $e^L = e$ and $f^L =
f$. For $r = 2$, consider $\mathcal{B}(2p \om) =
\tilde{\mathcal{B}}(2p \om)$:
$$u_{2p}\rightarrow u_{2(p-1)} \rightarrow u_{2(p-2)} \rightarrow
\cdots \rightarrow u_{2(2-p)}\rightarrow u_{2(1-p)}\rightarrow
u_{-2p}$$ which is decomposed in $\mathcal{B}^L(p\check\om)\sqcup
\mathcal{B}^L((p-1)\check\om)$ as a $^L\g$-crystal:
$$(u_{2p}\rightarrow u_{2(p-2)} \rightarrow \cdots \rightarrow
u_{-2p})\sqcup (u_{2p - 2}\rightarrow u_{2p - 6}\rightarrow \cdots
\rightarrow u_{2- 2p}).$$ 
Here Conjecture \ref{lchar} is just the elementary decomposition:
$$y^{2p} + y^{2p-2} + \cdots + y^{-2p} = (y^{2p} + y^{2p - 4} + \cdots
+ y^{-2p}) + (y^{2p - 2} + y^{2p - 6} + \cdots + y^{2 - 2p}).$$

We have the following Theorem, due to \cite{kmt2, kas1} (see for example
\cite[Theorem 2.1]{kas}):

\begin{thm}\label{reduction}
A finite $\Glie$-crystal $B$ is normal if and only if for any 
$J\subset I$ with at most two elements, $B$ is normal
as a $\Glie_J$-crystal.
\end{thm}

So it is of particular importance to study Lie algebras of rank
$2$. We will now prove Conjecture \ref{twiscrys} for Lie algebras of
rank $2$ (and $r \leq 2$). Such a Lie algebra will be denoted by
$(X,r_1,r_2)$, where $1 \leq r_1,r_2 \leq 2$ are the labels. We
consider all crystals $\mathcal{B}(\lambda)$ such that
$\tilde{\mathcal{B}}(\lambda)\neq \emptyset$. For $\Glie$ of type
$B_2$ it implies $\lambda\in P'$ but in general $\lambda$ is not
necessarily in $P'$.  For $(A_1\times A_1,2,2)$ and $(A_2,2,2)$ the
result is clear as we have $f_i^L = f_i$ and $e_i^L = e_i$. For types
$(A_1\times A_1,1,1)$, $(A_1\times A_1,1,2)$, $(A_1\times A_1,2,1)$,
the result follows from the rank $1$-case studied above. So we study
the two remaining case $(A_2, 1, 1)$ and $(B_2, 1, 2)$. In fact, we
prove

\begin{thm} Conjectures \ref{lchar} and \ref{twiscrys} hold for 
simply-laced $\Glie$ with $r = 2$ and for $B_2$.
\end{thm}

We cannot prove the statement for $B_2$ directly by using the result
for pseudo-fundamental representations (Proposition \ref{pfund}) as
the $e_i^L,f_i^L$ for the tensor product of ${}^L\g$-crystals do not
coincide with the operators defined from the tensor product of
$\g$-crystals.

\subsection{Type $(A_2, 1,1)$}\label{aidead}

Let $\lambda = R\om_1 + R'\om_2$ dominant in $P$. We have
$\lambda\equiv 0$, $\om_1$, $\om_2$ or $\om_1 + \om_2$ where $\equiv$
means mod $P'$ in this section.  Let $\Lambda = \{(i,1)|1\leq i\leq
R\}\cup\{(i,j)|R+1\leq i\leq R + R',1\leq j\leq 2\}$. Then
$\mathcal{B}(\lambda)$ is isomorphic \cite{kn} to the crystal of
tableaux $(T_{i,j})_{(i,j)\in \Lambda}$ with coefficients in $\{1,2,3\}$
which are semi-standard (i.e., $T_{i,j}\leq T_{i+1,j}$ and any $i,j$,
and $T_{i,1} > T_{i,2}$ for $i\geq 2R+1$).  Let
$$T_\lambda = \begin{pmatrix}
   &      & &1&\cdots &1
\\1&\cdots&1&2&\cdots &2
\end{pmatrix}$$ 
be the highest weight tableaux. Let us compute the tableaux
$T\in\tilde{\mathcal{B}}(\lambda)$ of highest weight for
$e_1^L = e_1^2$ and $e_2^L =
e_2^2$. $T = T(a,b,c)$ is characterized by $a$, $b$, $c$ such
that
$$T_{i,1} = \begin{cases} 1 \text{ for $i\leq a-1$,} \\2 \text{ for
$a\leq i\leq b-1$,} \\ 3 \text{ for $b\leq i$,}\end{cases}\text{ and }
T_{i,2} = \begin{cases} 1\text{ for $R+1\leq i\leq c-1$,} \\2\text{
for $c\geq i$.}
\end{cases}$$

\noindent The condition $e_1^2 T = 0$ is equivalent to the following:
$$(R+R' = c \text{ and }b - a \leq R' - 1)\text{ or }(R + R' =
c - 1\text{ and }b - a \leq R' + 1).$$
The condition $e_2^2 T = 0$ is equivalent to $c -
b\in\{0,1\}$. We have four cases:

1) $c = b = R + R'$ and $a\geq R + 1$. So $a = R + 1$. So
   $0 \equiv \mathrm{wt}(T) \equiv \lambda + \om_1 + \om_2$.

2) $c = b = R + R' + 1$ and $a\geq R$. If $a = R+1$ then $T =
   T_\lambda$ and $\lambda\equiv 0$. If $a = R$, then $\lambda \equiv \om_2$.
 
3) $c = b + 1 = R + R' $ and $a\geq R$. If $a = R+1$, then
   $\lambda \equiv \om_2$. If $a = R$, then $\lambda\equiv 0$.

4) $c = b + 1= R + R' + 1$ and $a\geq R - 1$. If $a = R+1$ then
   $\lambda \equiv \om_1$. If
   $a = R$ then $\lambda \equiv \om_1 + \om_2$. If $a = R - 1$ 
   then $\lambda \equiv \om_1$.

\noindent So for each value of $\lambda$ mod $P'$ we have $2$ highest
weight vectors $T_0$, $T_1$ of respective connected component
$\mathcal{B}'$ and $\mathcal{B}''$. We prove that $\mathcal{B}'\neq
\mathcal{B}''$ and that they are normal. This implies a stronger
result than Conjecture \ref{twiscrys}, that in this case
$\tilde{B}(\lambda)$ is normal.
$$(T_0,T_1) = \begin{cases} (T_\lambda, f_1f_2f_1f_2
T_\lambda) &\text{ if $\lambda \equiv 0$,}
\\(f_2T_\lambda ,f_2f_1f_1 T_\lambda) &\text{ if $\lambda \equiv \om_1$,}
\\(f_1 T_\lambda, f_1f_2f_2 T_\lambda)&\text{ if $\lambda \equiv \om_2$,}
\\(f_1f_2 T_\lambda, f_2f_1 T_\lambda)&\text{ if $\lambda \equiv \om_1 +
\om_2$.}
\end{cases}
$$

\noindent Note that $\mathrm{wt}(T(a,b,c))\in P'$ if and only if $b
\equiv 1 + R' [2]$ and $a\equiv c[2]$.

Let us treat in detail the case $\lambda \equiv 0$. We know by
Theorem \ref{dcrys} that $\mathcal{B}'$ is normal. In particular,
$\mathcal{B}'\neq \mathcal{B}''$. So we only have to prove that
$\mathcal{B}''(\lambda)$ is isomorphic as a $^L\g$-crystal to
$\mathcal{B}^L(\lambda')$ where $\lambda' = \Pi(\mathrm{wt}(T_1)) =
(R/2-1)\om_1 + (R'/2-1)\om_2$. We have $T_1 = T(R, R + R' -1, R + R')$
that is,
$$T_1 = \begin{pmatrix} & & & & 1 & \cdots & 1 & 1 & 2 \\ 1 &\cdots &
      1 & 2 & 2 & \cdots & 2 & 3 & 3\end{pmatrix}.$$ Let $\Lambda' =
      \{(i,1)|1\leq i\leq R/2-1\}\cup\{(i,j)|R/2+1\leq i\leq R/2 +
      R'/2-1,1\leq j\leq 2\}$. Then $\mathcal{B}^L(\lambda')$ is
      isomorphic to the crystal of semi-standard tableaux
      $(T_{i,j})_{(i,j)\in \Lambda'}$ with coefficients in
      $\{1,2,3\}$. For such a tableaux we define $a,b,c$ as
      above. Then consider
      $$\phi : T(a,b,c)\in\mathcal{B}^L(\lambda')\rightarrow T(2a,
      2b+1, 2c+2)\mathcal{B}''(\lambda).$$ Then $\phi$ is an
      isomorphism of $^L\g$-crystals.  First for $(R/2 , (R+R')/2
      -1, R/2 + R'/2-1)$ we get $T(R, R+R - 1, R + R') = T_1$.  Then
      it suffices to prove that $\phi(f_i T) =
      f_i^2\phi(T)$. Let $T = T(a,b,c)$.

For $f_1$ : if $R  + b \geq a + c$ and $a\geq 2$, then
$f_1 T = T(a - 1, b, c)$. We have $(R + 1) +
(2b + 1)\geq 2a + (2c + 2)$, so $f_1 \phi(T) = 
T(2a - 1, 2b+1, 2c+2)$. But we have also $(R + 1) + (2b + 1)\geq (2a -
1) + (2c + 2)$ so $f_1^2 \phi(T) = T(2(a - 1),
2b+1, 2c+2)$. 

If $R  + b \geq a + c$ and $a = 1$, then $f_1 T = 0$. We have
$(2R + 1) + (2b + 1)\geq 2a + (2c + 2)$, so $f_1 \phi(T) = 
T(1, 2b+1, 2c+2)$. But we have also $(R + 1) +(2b +
1)\geq (2a - 1) + (2c + 2)$ so $f_1^2 \phi(T) = 0$.

If $R + b < a + c$ and $c > b$ then $f_1 T$ corresponds to
$(a, b, c - 1)$. We have $(2R + 1) + (2b + 1)\leq 2a + 2c < 2a + 2c +
2$, so $f_1 \phi(T) = T(2a , 2b+1, 2c + 1)$. But
we have also $(R + 1) + (2b + 1) < 2a  + (2c + 1)$ so $f_1^2
\phi(T) = T(2a , 2b+1, 2(c-1)+2)$.

If $R + b < a + c$ and $c = b$ then $f_1 T = 0$. We have $(2R
+ 1) + (2b + 1)\leq 2a + 2c < 2a + 2c + 2$, so $f_1 \phi(T) 
= T(2a , 2b+1, 2c + 1)$. But we have also $(R + 1) + (2b
+ 1) < 2a  + (2c + 1)$ so $f_1^2 \phi(T) = 0$.

For $f_2$ : if $b > a$ then $f_2 T = T(a
, b - 1, c)$. We have $2b + 1\geq 2a + 3 > 2a$, so $f_1
\phi(T) = T(2a, 2b, 2c+2)$. But we have also $2b > 2a$,
so $f_1^2 \phi(T) = T(2a, 2b-1, 2c+2)$.

If $b = a$ then $f_2 T = 0$. We have $2b + 1 = 2a + 1 > 2a$,
so $f_1 \phi(T) = T(2a, 2b, 2c+2)$. But then $2b
= 2a$, so $f_1^2 \phi(T) = 0$.

For the cases $\lambda\equiv \om_1$ or $\lambda \equiv \om_1 + \om_2$
we give only the formulas of isomorphisms of $^L\g$-crystals as
above.

Let $\lambda\equiv \om_1$ (the case $\lambda\equiv \om_2$ is
symmetric). $T_0 = T(R+1,R + R',R + R' + 1)$, $T_1 = T(R-1,R + R',R +
R' + 1)$. Let $\phi' : \mathcal{B}((R+1)/2\om_1 + (R'/2 -
1)\om_2)\rightarrow \mathcal{B}'$
$$\phi' : T(a,b,c)\mapsto 
\begin{cases}T(R+1,2b - 1, 2c) &\text{ if $a = (R+3)/2$,}
\\T(2a - 1,2b - 1,2c - 1)&\text{ if $a < (R+3)/2$,}
\end{cases}$$ 
and
$$\phi'': T(a,b,c)\in \mathcal{B}((R-1)/2\om_1 + R'/2 \om_2)\mapsto
T(2a,2b + 1,2(c + 1))\in\mathcal{B}''.$$

For $\lambda\equiv \om_1 + \om_2$, $T_0 = T(R+1,R + R',R + R')$, $T_1
= T(R,R + R',R + R' + 1)$. Let
$$\phi' : T(a,b,c)\in\mathcal{B}((R-1)/2\om_1 + (R' -
1)/2\om_2)\mapsto T(2a,2b,2c)\in\mathcal{B}',$$
$$\phi'': T(a,b,c)\in\mathcal{B}((R-1)/2\om_1 + ((R' - 1)/2
\om_2)\mapsto T(2a-1,2b,2c+1)\in \mathcal{B}''.$$

\begin{rem} In the course of the proof we have found the following
Langlands duality branching rules (see the end of Section
\ref{charstate}) for irreducible representations of $(A_2, 1,1)$ and
the symmetric ones: ($\lambda_1,\lambda_2 > 0$)
$$\Pi(\chi(2\lambda_1\om_1 + 2\lambda_2\om_2)) = \chi^L(\lambda_1
\check\om_1 + \lambda_2\check\om_2) + \chi^L((\lambda_1 -
1)\check\om_1 + (\lambda_2 - 1)\check\om_2),$$ 

$$\Pi(\chi(2\lambda_1\om_1)) = \chi^L(\lambda_1 \check\om_1)\text{ ,
}\Pi(\chi((2\lambda_1 - 1)\om_1)) = \chi^L((\lambda_1 - 1)
\check\om_1),$$

$$\Pi(\chi((2\lambda_1 + 1)\om_1 + 2\lambda_2\om_2)) =
\chi^L((\lambda_1 + 1) \check\om_1 + (\lambda_2 - 1)\check\om_2) +
\chi^L((\lambda_1 - 1)\check\om_1 + \lambda_2\check\om_2),$$

$$\Pi(\chi((2\lambda_1 - 1)\om_1 + 2(\lambda_2 - 1)\om_2)) = 2
\chi^L((\lambda_1 - 1) \check\om_1 + (\lambda_2 - 1)\check\om_2).$$
\end{rem}

\subsection{Application to symmetric cases}\label{appsc}
Consider a simply-laced $\Glie$ with $r = 2$.

\begin{prop}\label{symdeux} For $\lambda\in P'$,
  $\tilde{\mathcal{B}}(\lambda)$ is normal.
\end{prop}

In particular, Conjectures \ref{lchar} and \ref{twiscrys} hold for
these types. In fact, we have proved a stronger result as the crystal
is normal.

\medskip

\demo By Theorem \ref{reduction}, it suffices to prove the result for
the subalgebras of rank $2$.  For subalgebras of type $A_2$, the
statement follows from Section \ref{aidead}. For the subalgebras of
type $A_1\times A_1$, it suffices to prove that if $C_{i,j} = 0$ then
the $f_i^L,f_j^L$ commute. But it is clear as the
$f_i,f_j$ commute.  \qed

\subsection{Type $(B_2,1,2)$} 

Let $\lambda = 2R\om_1 + R'\om_2\in P'$.  Let $\Lambda = \{(i,1)|1\leq
i\leq 2R\}\cup\{(i,j)|2R+1\leq i\leq 2R + R',1\leq j\leq 2\}$. Then
$\mathcal{B}(\lambda)$ is isomorphic \cite{kn} to the crystal of
tableaux $(T_{i,j})_{(i,j)\in \Lambda}$ with coefficients in
$\{1,2,\overline{2},\overline{1}\}$ which are semi-standard (i.e.,
$T_{i,j}\preceq T_{i+1,j}$ and any $i,j$, and $T_{i,1} \succ T_{i,2}$
for $i\geq 2R+1$ for the ordering $1\preceq 2\preceq
\overline{2}\preceq \overline{1}$) and such that for $i\geq 2R+1$,
$(T_{i,1},T_{i,2})\neq (\overline{1},1)$ and $(T_{i+1,1},T_{i,2})\neq
(\overline{2}, 2)$.

Let $T_\lambda$ be the highest weight tableaux. The tableaux $T =
T_\epsilon(a,b,c,d)$ is characterized by $a$, $b$, $c$, $d$ and
$\epsilon\in\{0,1\}$ such that
$$T_{i,1} = \begin{cases}
1 \text{ for $i\leq a-1$,}
\\2 \text{ for $a\leq i\leq b-1$,}
\\ \overline{2} \text{ for $b\leq i\leq c-1$,}
\\ \overline{1} \text{ for $c\leq i$,}\end{cases}\text{ and }
T_{i,2} = \begin{cases}
1\text{ for $2R+1\leq i\leq c-\epsilon-1$,}
\\2\text{ for $c-\epsilon \geq i\leq d-1$,}
\\\overline{2}\text{ for $d\leq i$.}
\end{cases}$$
In fact, $\begin{pmatrix}2\\\overline{2}\end{pmatrix}$ appear at most
once (it can appear in $T_1$ and does not appear in $T_0$).

Let us compute the tableaux $T$ of highest weight for the operators
$e_1^L = e_1^2$ and $e_2^L = e_2$.
The condition $e_2T = 0$ implies $d = R + R' + 1$.  The
condition $e_1^2 T = 0$ implies $c = d = R + R' + 1$.  Let us
consider the $3$ classes of such tableaux:

Tableaux (A) : $T_{R,1} = \overline{1}$ (that is $c\leq
R$). $e_2 T = 0$ gives $R'\geq c - b$. $e_1^2 T = 0$
gives $R' = 0 = c - b$ and $2R\leq a$. So all coefficients are equal
to $1$ except $T_{2R,1}\in\{1,2,\overline{1}\}$.

Tableaux (B) : $T_{R,1} = \overline{2}$ (that is $c > R$ and $b\leq
R$). $e_2 T = 0$ gives $b = R + R' + 1 - \epsilon$. So
$\epsilon = 1$ and $R' = 1$. $e_1^2 T = 0$ gives $a \geq
2R$. So $(T_{R+1,1},T_{R+1,2}) = (\overline{2},2)$, $T_{R,1}\in\{1,
2\}$ and all other coefficients are equal to $1$.

Tableaux (C) : $T_{1,1} \preceq 2$ (that is $b > R$). $e_2 T =
0$ gives $b = R + R' + 1 - \epsilon$. Then $e_1^2 T = 0$ gives
$a\geq 2R$.

For $R' = 0$ and $R > 0$ : we get $3$ tableaux $T_\lambda$,
$f_1 T_\lambda$, $f_1f_2
f_1T_\lambda$.

For $R = 0$ and $R' > 0$ : we get $2$ tableaux $T_\lambda$,
$f_1f_2 f_1T_\lambda$.

For $R,R' > 0$ : we get $4$ tableaux 
$$T_\lambda = \begin{pmatrix}&&&1&\cdots&1
\\1&\cdots&1&2&\cdots&2
\end{pmatrix},$$ 
$$T_1 = f_1 T_\lambda = 
\begin{pmatrix}&&&&1&\cdots&1
\\1&\cdots&1&2&2&\cdots&2
\end{pmatrix},$$ 
$$T_2 = f_1f_2 T_\lambda = 
\begin{pmatrix}&&&1&\cdots&1&2
\\1&\cdots&1&2&\cdots&2&\overline{2}
\end{pmatrix},$$ 
$$T_3 = f_1f_2f_1 T_\lambda
= \begin{pmatrix}&&&&1&\cdots&1&2
\\1&\cdots&1&2&2&\cdots&2&\overline{2}
\end{pmatrix}.$$
We concentrate on the case $R,R' > 0$ (the cases $R = 0$ or $R' = 0$
can be easily deduced from it).  By Theorem \ref{dcrys} the connected
component of $T_\lambda$ is isomorphic to the crystal of a simple
$^L\Glie$-module. In particular it contains $T_1, T_2, T_3$.  Let
$\mathcal{B}$ (resp. $\mathcal{B}'$) be the union of the component of
$T_1, T_2$ (resp. the component of $T_3$).  We have $u\in \mathcal{B}$
if and only if $\mathrm{wt}(u)\in \lambda - (1 + 2\ZZ)\alpha_1 -
\ZZ\alpha_2$.  So the component $\mathcal{B}\cap \mathcal{B}' =
\emptyset$. In the monomial model
$\mathcal{M}(Y_{2,0}^{R'}Y_{1,1}^{2R})$, $T_3$ corresponds to
$Y_{1,1}(Y_{2,0}^{R'}Y_{1,1}^{2(R-1)})Y_{1,5}^{-1}$. By Theorem
\ref{ld} the $^L\Glie$-crystal generated by
$Y_{2,0}^{R'}Y_{1,1}^{2(R-1)}$ is the crystal of the simple
$^L\Glie$-module of highest weight $(R-1)\check\om_1 + R'
\check\omega_2$. Bu the multiplication by $Y_{1,1}Y_{1,5}^{-1}$ does
not change the action of the crystal operators here, and so
$\mathcal{B}_3$ is also isomorphic to this crystal.

For $\mathcal{B}$ we write explicitly the bijection by using the three
cases as above. To do it we also use the dual tableaux realization of
$\mathcal{B}^L (\mu)$ for $\mu = \mu_1 \check\om_1 +
\mu_2\check\om_2$.

Let $\Lambda^L = \{(i,2)|\mu_1 < i\leq \mu_1 + \mu_2\}\cup\{(i,j)|1\leq
i\leq \mu_1,1\leq j\leq 2\}$. $\mathcal{B}^L(\mu)$ is isomorphic
\cite{kn} to the crystal of tableaux $(T_{i,j})_{(i,j)\in \Lambda^L}$ with
coefficients in $\{1,2,\overline{2},\overline{1}\}$ which are
semi-standard as above.  The tableaux $T = T_\epsilon^l(a,b,c,d)$ is
characterized by $a$, $b$, $c$, $d$ and $\epsilon\in\{0,1\}$ such that
$$T_{i,1} = \begin{cases}
2 \text{ for $i\leq a-1$,}
\\\overline{2} \text{ for $a\leq i\leq b-1$,}
\\ \overline{1} \text{ for $b\leq i\leq \mu_1$,}
\end{cases}\text{ and }
T_{i,2} = \begin{cases}
1\text{ for $i\leq b - \epsilon - 1$,}
\\2\text{ for $b - \epsilon \geq i\leq c - 1$,}
\\\overline{2}\text{ for $c\leq i\leq d -1$,}
\\\overline{1}\text{ for $d\leq i$.}
\end{cases}$$
Let $\mathcal{B}_1^L = \mathcal{B}^L(R\check\om_1 + (R' -
1)\check\om_2)$ and $\mathcal{B}_2^L = \mathcal{B}^L((R -
1)\check\om_1 + (R'+ 1)\check\om_2)$. We define $\Psi :
\mathcal{B}_1^L\sqcup \mathcal{B}_2^L\rightarrow \mathcal{B}$.  The
general idea to define the map is to replace $(1,1)$, $(2,2)$,
$(2,\overline{2})$, $(\overline{2},\overline{2})$,
$(\overline{1},\overline{1})$ in the first part of the tableaux
respectively by $\begin{pmatrix}1\\2\end{pmatrix}$,
$\begin{pmatrix}1\\\overline{2}\end{pmatrix}$,
$\begin{pmatrix}2\\\overline{2}\end{pmatrix}$,
$\begin{pmatrix}2\\\overline{1}\end{pmatrix}$,
$\begin{pmatrix}\overline{2}\\\overline{1}\end{pmatrix}$, and to
replace $\begin{pmatrix}1\\2\end{pmatrix}$,
$\begin{pmatrix}1\\\overline{2}\end{pmatrix}$,
$\begin{pmatrix}2\\\overline{2}\end{pmatrix}$,
$\begin{pmatrix}2\\\overline{1}\end{pmatrix}$,
$\begin{pmatrix}\overline{2}\\\overline{1}\end{pmatrix}$ in the second
part of the tableaux respectively by $(1)$, $(2)$, $()$,
$(\overline{2})$, $(\overline{1})$.  In general, it cannot be done in
the obvious way as other term may appear as $(1,2)$,
$(2,\overline{1})$ and so we have to do the following case by case
description.

Tableaux (C). $T_\pm(a,b,c,d)\in \mathcal{B}\Leftrightarrow a \equiv
\epsilon [2]$. Let $\beta\geq R + 1$. We set :
$$T_0(\alpha,\beta,\gamma,\delta)\in\mathcal{B}_1^L\mapsto T_1(2\alpha
- 1, \beta + R, 1 + R + \gamma, 1 + R + \delta),$$
$$T_0(\alpha,\beta,\gamma,\delta)\in\mathcal{B}_2^L\mapsto T_0(2\alpha
, \beta + R , R + \gamma, R + \delta).$$

Tableaux (B). $T_\pm(a,b,c,d)\in \mathcal{B}\Leftrightarrow a \equiv
\epsilon [2]$. Let $\beta \leq R < \gamma$. We set
$$T_\epsilon(\alpha,\beta,\gamma,\delta)\in\mathcal{B}_1^L\mapsto
T_1(2\alpha - 1, 2\beta - 1 - \epsilon, 1 + R + \gamma, 1 + R +
\delta),$$
$$T_\epsilon(\alpha,\beta,\gamma,\delta)\in\mathcal{B}_2^L\mapsto
T_0(2\alpha , 2 \beta - \epsilon , R + \gamma, R + \delta).$$

Tableaux (A). $T_0(a,b,c,d)\in \mathcal{B}\Leftrightarrow c\equiv a +
1[2]$. Let $\gamma\leq R$. We set :
$$T_\epsilon(\alpha,\beta,\gamma,\delta)\in\mathcal{B}_1^L\mapsto
T_0(2\alpha - 1, 2\beta - \epsilon - 1 , 2\gamma, 1 + R + \delta),$$
$$ T_\epsilon(\alpha,\beta,\gamma,\delta)\in\mathcal{B}_2^L\mapsto
\begin{cases}
T_0(2\alpha, 2 \beta - \epsilon , 2\gamma - 1, R + \delta)&\text{ if
($\epsilon = 1$ or $\beta < \gamma$) and $\delta > R$}, \\T_0(2\alpha
- 1, 2 \beta, 2\beta, R + \delta)&\text{ if $\epsilon = 0$, $\beta =
\gamma$ and $\delta > R$,} \\T_0(2\alpha - 1, 2 \beta - 1 - \epsilon,
2\gamma, 2R + 1)&\text{ if $\delta = R$.}
\end{cases}
$$
It is straight forward to check that the properties of Conjecture
\ref{twiscrys} are satisfied.

\begin{rem} In the course of the proof we have found the following
Langlands duality branching rules for
irreducible representations of $(B_2, 2,1)$: ($\lambda_1,\lambda_2 >
0$)
$$\Pi(\chi(2\lambda_1\om_1 + \lambda_2\om_2)) = \chi^L(\lambda_1
\check\om_1 + \lambda_2\check\om_2) + \chi^L(\lambda_1\check\om_1 +
(\lambda_2 - 1)\check\om_2)$$
$$+ \chi^L((\lambda_1 - 1)\check\om_1 + (\lambda_2
+ 1)\check\om_2) + \chi^L((\lambda_1 - 1)\check\om_1 +
\lambda_2\check\om_2),$$ 
$$\Pi(\chi(2\lambda_1\om_1)) = \chi^L(\lambda_1
\check\om_1)  + \chi^L((\lambda_1 - 1)\check\om_1 + \check\om_2) +
\chi^L((\lambda_1 - 1)\check\om_1),$$

$$\Pi(\chi(\lambda_2\om_2)) = \chi^L(\lambda_2\check\om_2) +
\chi^L((\lambda_2 - 1)\check\om_2).$$
\end{rem}

\subsection{A proposed deformation process} Suppose that $r = 2$. 
We have proved the statement of Conjecture \ref{twiscrys} for rank
$2$, but we cannot use Theorem \ref{reduction} directly for general
rank. For example, for type $B_3$, $\tilde{\mathcal{B}}(\lambda)$ is a
normal crystal for ${}^L\g_{\{1,2\}}$ and ${}^L\g_{\{1,3\}}$. We use
the rank $2$ to deform the $3$-arrows so that we get
${}^L\g_{\{2,3\}}$. But then we may not preserve the
${}^L\g_{\{1,3\}}$-crystal structure.

We propose a conjectural inductive process to redefine the crystal
operators of $\tilde{\mathcal{B}}(\lambda)$ so that we get a normal
crystal. Suppose that we know the result for rank lower than $n - 1$
for an $n\geq 3$.  Let $I = I_1\sqcup I_2$ where $I_k = \{i\in I|r_i =
k\}$.  We assume $|I_2|\geq 2$ (the case $|I_1|\geq 2$ can be treated
in a symmetric way by Proposition \ref{symdeux}).  We use the notation
$I_2 = \{1,\cdots, i_0\}$ and $I_1 = \{i_0 + 1,\cdots, n\}$ so that
$C_{i_0,i_0 + 1} = -1$. Let $\overline{I_2} = I_2-\{i_0\}$,
$\overline{I_1} = I_1 \cup\{i_0\}$.

Let $\lambda\in P'$ and fix a class $P'' = \mu + Q''\subset \lambda +
Q$ mod $Q'$.  Then $\mathcal{B} =
\{u\in\tilde{\mathcal{B}}(\lambda)|\mathrm{wt}(u)\in P''\}$ is a union
of connected component of $\tilde{\mathcal{B}}(\lambda)$ as the weight
of the vectors in a connected component are in the same class.  For
$\mu_1, \mu_2\in P''$, we have $\mu_1 - \mu_2 = \sum_{i\in I}n_i
\alpha_i^L$ where $n_i\in\ZZ$ and the $\alpha_i^L$ are the simple
roots of $^L\g$. We put $N(\mu_1,\mu_2) = \sum_{i\in I}n_i$.  Let
$\mu'\in\{\mathrm{wt}(u)|u\in\mathcal{B}\}$ such that $N(\mu,\mu')$ is
maximal. It is well defined, that is to say independent of the choice
of $\mu$, as for $\mu_1,\mu_2,\mu_3\in P''$ we have $N(\mu_1,\mu_2) +
N(\mu_2,\mu_3) = N(\mu_1,\mu_3)$. We set $N(\mu_1) = N(\mu',\mu_1)$.
For $N\geq 0$, let $W_N = \{u\in\mathcal{B}|N(\mathrm{wt}(u)) = N\}$.

For $\mathcal{C}$ a (normal) $^L\g$ crystal, by truncated (normal) crystal of
$\mathcal{C}$ we mean for a certain $N\in\ZZ$ the set
$\{u\in\mathcal{C}|N(\mathrm{wt}(u))\geq N\}$ with the maps
$\mathrm{wt}^L, e_i^L,\epsilon_i,\phi_i$ restricted to it and
the map $f_i^L$ restricted to
$\{u\in\mathcal{C}|N(\mathrm{wt}(u))\geq N-1\}$.

To start we set all $(f_i^L)' = f_i^L$, $(e_i^L)' = e_i^L$.  By
induction on $N\geq 0$, we redefine $(f_i^L)'$ on $\sqcup_{M\leq N -
1}W_M$ (or equivalently $(e_i^L)'$ on $\sqcup_{M\leq N}W_M$). We say
that the process does not fail if $(\sqcup_{M\leq N}W_M,
\mathrm{wt}^L, \epsilon_i^L, \phi_i^L, (f_i^L)', (e_i^L)')$ is a
normal truncated crystal.

For $N = 0$ we do not change the maps. Let $N\geq 0$.

Let $i\in I_1$ and $u\in W_{N-1}$ such that $\exists j\in
\overline{I_2}, \epsilon_j^L(u) > 0$.  If $\phi_i^L(u) = 0$ we set
$(f_i^L)'(u) = f_i^L(u) = 0$. Otherwise let $v = e_j^L(u)\neq 0$. Then
$\phi_i^L(v) = \phi_i^L(u)\neq 0$ so $w = (f_i^L)'(v)\neq 0$. Then
$\phi_j^L(w) = \phi_j^L(v)\neq 0$ so $x = f_j^L(w)\neq 0$. We set
$(f_i^L)'(u) = x$.
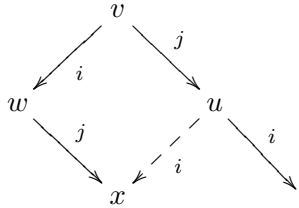
$$\xymatrix{& v \ar[dl]^i\ar[dr]^j & & \\ w \ar[dr]^j& &
u\ar@{-->}[dl]^i\ar[dr]^i& \\ &x&&}$$ 
We have
$\epsilon_i^L(x) = \epsilon_i^L(w) = \epsilon_i^L(v) + 1 =
\epsilon_i^L(u) + 1$ and $(f_i^L)'(u)$ is well-defined
(independent on $j\in \overline{I_2}$).
$(e_i^L)'(y)$ is now defined for 
$y\in W_N$ such that $\sum_{j\in \overline{I}_2}\epsilon_j^L(u) > 0$.

Let $\mu\in P'$ and $\mathcal{U}_\pm = \{y\in (W_N)_\mu|\pm\sum_{j\in
\overline{I_2}} \epsilon_j^L(y) \leq \pm 1/2\}$. We redefine
$(e_i^L)'$ on $\mathcal{U}_+$ by induction on $i\geq i_0$.  Let $u\in
\mathcal{B}_i = \{u\in (W_N)_{\mu + \alpha_i}|\phi_i^L(u) > 0 \text{ ,
}u\notin (e_i^L)'(\mathcal{U}_-)\}$.  Consider the truncated
$^L\g$-crystal $\sqcup_{M\leq N-1}W_M$ and $\mathcal{C}$ be the
corresponding normal crystal with the injection $\Psi : \sqcup_{M\leq
N-1}W_M\rightarrow \mathcal{C}$. We have $\phi_i^L(\Psi(u)) =
\phi_i^L(u) > 0$ and so $v = f_i^L(\Psi(u))\neq 0$.  If there is
$i_0\leq j\leq i - 1$ such that $\epsilon_j^L(v) > 0$, let $w =
e_j^L(v)$. We set $(f_i^L)'(u) = (f_j^L)'\Psi^{-1}(w) = x$.
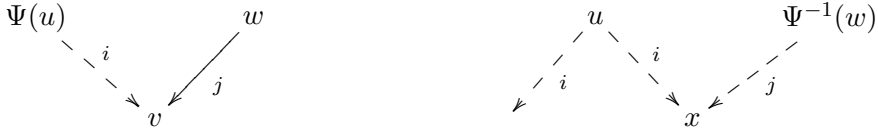
$$\xymatrix{\Psi(u) \ar@{-->}[dr]^i & & w\ar[dl]^j&&
&&u\ar@{-->}[dl]^i\ar@{-->}[dr]^i&&\Psi^{-1}(w)\ar@{-->}[dl]^j \\ & v &&&
&& & x & }
$$
As $|\overline{I_1}| < n$ we have : 
$$|\{v\in\mathcal{B}_i|(\epsilon_j^L(f_i^L(\Psi(u))))_{j\in
\overline{I_1}} = (a_j)_{j\in \overline{I_1}}\}| = |\{v\in
\mathcal{U}_+ |(\epsilon_j^L(v))_{j\in \overline{I_1}} = (a_j)_{j\in
\overline{I_1}}\}|$$ for a given $(a_j)_{j\in \overline{I_1}}$. So we
can define $(f_i^L)'(u)$ for $u\in\mathcal{B}_i$ such that
$\sum_{i_0\leq j\leq i-1}\epsilon_j(v) = 0$. We get $(f_i^L)' :
\mathcal{B}_{i_0} \rightarrow \mathcal{U}_+$ injection. Moreover
conjecturally for $i = i_0$, we can choose $(f_{i_0}^L)'$ compatible
with the $\phi_i^L$, $i\in \overline{I_2}$ (in other words, there is
"enough dimension" in weight spaces to do it) and then we can redefine
$(f_{i_0}^L)'$ on $\mathcal{U}_+$ so that the structure of
${}^L\Glie_{I_2}$-crystal is not modified.

If the conjectural point is satisfied, the process never fails, and
the new crystal is normal for any $^L\g_J$ where $|J|\leq 2$. Then we
can conclude with Theorem \ref{reduction}.

\end{document}